\documentclass[pdf,aos,preprint]{imsart}
\usepackage{xr}
\usepackage{amsmath, mathrsfs}
\RequirePackage[numbers]{natbib}
\usepackage{graphicx}
\usepackage{color} 
\usepackage{verbatim}
\usepackage{hyperref}
\newcommand{\qed}{{\unskip\nobreak\hfil\penalty50\hskip2em\vadjust{}
            \nobreak\hfil$\Box$\parfillskip=0pt\finalhyphendemerits=0\par}}

\newtheorem{thm}{Theorem}[section] 
\newtheorem{lemma}{Lemma}[section] 
\newtheorem{cor}{Corollary}[section]

\newtheorem{definition}{Definition}[section]
\newcommand{\bed}{\begin{definition}}
\newcommand{\eed}{\end{definition}}
\newcommand{\beq}{\begin{equation}}
\newcommand{\eeq}{\end{equation}}

\newcommand{\eps}{\epsilon}

\newcommand{\bitem}{\begin{itemize}}
\newcommand{\eitem}{\end{itemize}}

\newcommand{\goto}{\rightarrow}

\newcommand{\mexp}{\mathrm{exp}}

\newcommand{\beqn}{\begin{equation}}
\newcommand{\eeqn}{\end{equation}}
\newcommand{\balign}{\begin{align}}
\newcommand{\ealign}{\end{align}}

\newcommand{\diag}{\mathrm{diag}}

\newcommand{\bel}{\begin{eqnarray}\label}
\newcommand{\eel}{\end{eqnarray}}
\newcommand{\bes}{\begin{eqnarray*}}
\newcommand{\ees}{\end{eqnarray*}}

\newcommand{\tr}{\mathrm{tr}}

\usepackage{amssymb}
\begin{document}
\begin{frontmatter}

\title{Detecting Rare and Weak spikes in large covariance matrices}
\begin{aug}
\author{\fnms{Zheng Tracy} \snm{Ke}\thanksref{t1}\ead[label=e1]{zke@fas.harvard.edu}}

\affiliation{
Harvard University\thanksmark{t1}
} 
	
	\address{Z. Ke\\
		Department of Statistics\\
		Harvard University\\
		Cambridge, Massachusetts, 02138\\
		USA\\
		\printead{e1}}
\end{aug}

\begin{abstract} 
Given $p$-dimensional Gaussian vectors  $X_i \stackrel{iid}{\sim} N(0, \Sigma)$, $1 \leq i \leq n$, where $p \geq n$,    we are interested in testing 
a null hypothesis where $\Sigma = I_p$ against an alternative hypothesis where all eigenvalues of $\Sigma$ are $1$, except for $r$ of them are larger than $1$  (i.e.,  spiked eigenvalues).   
 
We consider a Rare/Weak setting where the spikes are sparse (i.e., $1 \ll r \ll p$) and  individually weak  (i.e., each spiked eigenvalue is only slightly larger than $1$), and discover a  phase transition:  the two-dimensional phase space that calibrates the spike sparsity and strengths partitions into the {\it Region of Impossibility} and the {\it Region of Possibility}. 
In Region of Impossibility, all tests are (asymptotically) powerless in separating the alternative from the null.   
In  Region of Possibility, there are tests that have (asymptotically) full power. 

We consider a CuSum test, a trace-based test, an eigenvalue-based Higher Criticism test,  and a Tracy-Widom test \cite{Johnstone2001}, and show that the first two tests have asymptotically full power in Region of Possibility.

To use our results from a different angle,  we 
derive new bounds for (a) empirical eigenvalues,  and (b)  cumulative sums of the empirical eigenvalues, both under the alternative hypothesis. 
Part (a) is related to those in  \cite{BAP, Paul2007},  but both the settings and results are different.

The study requires careful analysis of the $L^1$-distance of our testing problem and  delicate Radom Matrix Theory.  
Our technical devises include (a)  a Gaussian proxy model, (b) Le Cam's  comparison of experiments, and (c) large deviation bounds on empirical eigenvalues. 

\end{abstract}

 

\end{frontmatter}

\section{Introduction} 
\label{sec:intro} 
\setcounter{equation}{0}
Suppose we have $n$ normal vectors   $X_i \in \mathbb{R}^p$: 
\begin{equation} \label{Gaussianmodel}  
X_i \stackrel{iid}{\sim} N(0, \Sigma), \qquad \Sigma \in \mathbb{R}^{p,p}, \qquad  1 \leq i \leq n.  
\end{equation}  
We are interested in testing the null hypothesis that $\Sigma$ is the $p$ by $p$  identity matrix $I_p$ against an alternative hypothesis that  $\Sigma$ is a so-called spike matrix \cite{Johnstone2001}: for an integer $r$ such that $1 \leq r \ll p$,  all eigenvalues of $\Sigma$ are $1$, except for the top $r$ of them are  larger than $1$ (i.e., each of them is  spiked): 
\begin{equation} \label{latenroots} 
\ell_1 \geq \ell_2\geq \ldots \geq \ell_r > 1. 
\end{equation} 
Motivated by the recent interest of ``$p > n$",   we assume 
$p \geq n$,   
but this should not be  taken as a constraint. Let $X \in \mathbb{R}^{n,p}$ be the data matrix so that 
$X' = [X_1, X_2, \ldots, X_n]$, and denote the empirical covariance matrix by 
\[
\hat{\Sigma} = (1/n) X'X. 
\]
Since $p \geq n$, with probability $1$ \cite{Uhlig94},  $\hat{\Sigma}$ has $n$ distinct positive eigenvalues 
\[ 
 \lambda_1 >   \lambda_2  >   \ldots >   \lambda_n > 0.   
\]

This testing problem is of interest in many application areas. 
\begin{itemize}
\item {\it Covert communication}. In computer security and privacy, covert channels are widely used.   Consider a channel with $r$ antennas at transmitter and $p$ antennas at receiver. The output vectors satisfy $X_i= Hy_i + Z_i$, $1\leq i\leq n$, where $y_i\overset{iid}{\sim} N(0, \delta \cdot I_r)$ are the input vectors, $Z_i\overset{iid}{\sim} N(0, I_p)$ is noise, and $H\in\mathbb{R}^{p,r}$ is a confidential ``channel matrix" \cite{Verdu}. When using a covert channel, we would like to know whether our ``enemy" can notice that we are sending signals. Note that $X_i$'s follow Model \eqref{Gaussianmodel} with $\Sigma=I_p +\delta HH'$. From the perspective of our enemy, $H$ is unavailable, so it becomes a problem of detecting spiked eigenvalues in $\Sigma$. We will revisit this application in Section~\ref{subsec:covert}.
\item {\it Inference of genetic population structure.} One of main challenges in analyzing genetic data is to explore whether the samples contain subpopulations that are genetically distinct \cite{price2006principal,patterson2006population}. Let $Y\in\{0,1,2\}^{n,p}$ be the matrix of SNP counts for $p$ markers and $n$ subjects. When there are no sub-populations, $\{Y(i,j)\}_{i=1}^n$ are modeled as $iid$ samples from $Binomial(2,p_j)$, where $p_j\in (0,1)$ is the Minor Allele Frequency of marker $j$. Consider the normalized data matrix $X\in\mathbb{R}^{n,p}$, where $X(i,j)=[Y(i,j)-2\hat{p}_j]/\sqrt{2\hat{p}_j(1-\hat{p}_j)}$ and $\hat{p}_j=(2n)^{-1}\sum_{i=1}^n Y(i,j)$. Then, entries of $X$ are independent, with (approximately) zero mean and unit variance. When there are subpopulations, rows of $X$ are confounded by latent ancestry variables, and its covariance matrix has spiked eigenvalues \cite{patterson2006population}. Inference of genetic sub-populations reduces to detecting spikes in the covariance matrix. 
\item {\it Multiple testing.} How to account for data dependence is a central challenge in large-scale multiple testing \cite{efron2007correlation,DJ15}. Recent works \cite{leek2007capturing, friguet2009factor, fan2012estimating} used Model \eqref{Gaussianmodel}-\eqref{latenroots} to model test statistics, where each spiked eigenvalue comes from an unobserved confounding variable. Under this framework, various factor-adjusted multiple testing procedures have been proposed. Then, a fundamental problem is to detect the existence of confounding factors, so that we know when to use factor-adjusted methods instead of classical ones. 
\item {\it Network community detection.} Given a large social network, we are often interested in testing whether it contains only one community or multiple communities. Let $A\in\{0,1\}^{n,n}$ be the adjacency matrix of a symmetrical network with $n$ nodes. Consider the centered and re-scaled adjacency matrix $\widetilde{A}\in\mathbb{R}^{n,n}$ by $\widetilde{A}(i,j)=[A(i,j)-\hat{p}]/\sqrt{\hat{p}(1-\hat{p})}$, where $\hat{p}=\sum_{i<j}A(i,j)/[n(n-1)/2]$. Under the stochastic block model, it is well-known that $\widetilde{A}$ is approximately a Wigner matrix when there is only one community and has a few spiked eigenvalues when there are multiple communities \cite{bickel2016hypothesis, lei2016goodness}. 
\end{itemize}
Other applications include low-rank matrix recovery \cite{CandesPlan},  PCA and sparse PCA \cite{ZHT, Paul2007}, and high-dimensional clustering \cite{JKW, IFPCA}. 

In Model~\eqref{Gaussianmodel}-\eqref{latenroots}, we assume all eigenvalues of $\Sigma$ are $1$ under the null. This is a mild assumption since the data matrix is always pre-normalized in above applications. The setting that $\hat{\Sigma}$ is a Gaussian covariance matrix is not exactly the same as those in genetics (where data are sub-Gaussian) and in social networks (where we deal with Wigner matrices). But due to eigenvalue universality, the asymptotic behavior of many eigenvalue-based tests remains the same. Our testing framework can be viewed as an idealization of these applications, but it still captures the essential features.

Recently, encouraging progresses have been made to understanding this problem, such as Johnstone \cite{Johnstone2001}, Onatski et al. \cite{Hallin2013, Hallin2014} and Johnstone and Nadler \cite{Johnstone2015}.  
However, these works have been largely focused on the case where the number of spikes $r$ is small. In fact, 
in the asymptotic framework used in these papers,   $r$ 
is   fixed as $n \goto \infty$.  

In this paper, the primary interest is to consider the testing problem in the case where we have {\it many  weak  spikes}. 
For covert communication, this means the transmitter has a number of antennas, each with a relatively small capacity. For genetical data, it is the case where the subjects come from a number of sub-populations whose mutual genetic distinction is so weak that to separate them is subtle or even impossible; this can happen for studies conducted on European populations (say). Similar scenarios also happen in social networks: a large network often contains many ``weak communities" \cite{radicchi2004defining} that are hard to identify. Motivated by these applications,   
we adopt a Rare and Weak Spike (RWS) model where as $n \goto \infty$,  
\[
\mbox{$r \goto \infty$ in an algebraic rate, each spiked eigenvalue is   {\it slightly} larger than $1$}. 
\] 
The main contribution of the paper is three-fold.  
\begin{itemize} 
\item ({\it Phase transition}). We  discover that the 
two-dimensional phase space calibrating the spike sparsity and spike strengths partitions into 
the Region of Impossibility and the Region of Possibility. In the former, the spikes are so rare/weak that it is impossible to separate the alternative from the null. In the latter, the spikes are sufficiently strong and it is possible to separate the  alternative from the null.  
\item ({\it Methods}).  We propose a CuSum test and an (eigenvalue-based) Higher Criticism (HC) 
test\footnote{The HC test is based on $\lambda_i$ and is very different from those in literature.} as new approaches to the testing problem.  In particular, we show that the CuSum test is optimal.  
\item ({\it Bounds on eigenvalues}).  Knowledge/innovation flow is rarely a one-way street:  
while Random Matrix Theory (RMT) helps us  establish the lower/upper bounds of our testing problem, the latter also shed lights on some problems in RMT in return:        
we obtain new bounds on $\lambda_i$ and on the cumulative sums of $\lambda_i$,   a topic of great interest in RMT.  
\end{itemize} 
Our study requires delicate analysis, and the following are some noteworthy points of our technical contributions. 
\begin{itemize} 
\item The RWS model is similar to that in \cite{Johnstone2015, Hallin2014} and is known to be hard to analyze. We overcome the technical hurdle by using a Gaussian proxy model, which is comparably easier to analyze.  Our result on the intimacy of two models is of interest for its own sake, and can be used as a technical device to study problems other than spike detection.  
\item Le Cam's comparison of experiments provides a useful tool for analyzing complicated models. Using Le Cam's idea,  
we extend our main results on the RWS model to more complicated settings. 
\item Our analysis needs many recent results in RMT.  These results scatter across the literature 
and are presented in forms that are not always easy to access.  
With substantial efforts, we  adapt such results to our testing problem and make them more  accessible for us (see Section~\ref{sec:RMT}).  
\end{itemize} 

\subsection{Four test statistics}
\label{subsec:fivetests}  
We propose CuSum (CS) and (eigenvalue-based) Higher Criticism (HC) as two new tests. We also investigate a trace-based test and a Tracy-Widom (TW) test.

The trace-based test uses the empirical moments of $\lambda_i$ for testing. While there are many tests  of this kind, we focus our study on the following test: 
\[
S_n   =   \sum_{i = 1}^n \lambda_i   =  \tr(\hat{\Sigma}).\footnote{For any square matrix $A$,  $\tr(A)$ stands for the trace.} 
\] 
In Section \ref{subsec:main} we discuss other versions  of trace-based statistics.

The trace statistic uses the total sum of all eigenvalues.  A natural alternative is to use the  cumulative sum (CuSum) 
\[
S_k =  \sum_{i = 1}^k \lambda_i,  \qquad 1 \leq k \leq n,    
\]
where we hope for some $k < n$, $S_k$ contains stronger evidence against the null than $S_n$.   
Of course, we don't know how to pick the ``best" $k$.  To address this problem, we propose the following test statistic which we call the {\it CuSum}:  
\[
CS_n^* = \max_{1 \leq k \leq n} \{CS_{n,k}\}, \qquad \mbox{where $CS_{n,k} = \frac{S_k - E_0[S_k]}{SD_0(S_k)}$}.  
\]
Here, $E_0[S_k]$ and $SD_0[S_k]$ are the mean and standard deviation of $S_k$ under the null, respectively.  Such quantities do not have  a closed-form expression, but  can be conveniently simulated. 
In some cases, especially when the spikes are very rare, CuSum improves the trace-based test; see Section \ref{sec:Simul}.


In a similar spirit, we have the eigenvalue-based {\it Higher Criticism (HC)} test.\footnote{The test was briefly mentioned in a survey paper by Donoho and Jin \cite{DJ15}, but has not yet been formally studied.}  The HC statistic is defined by 
\begin{equation} \label{DefineHC} 
HC_n^* =  \max_{1 \leq k \leq n} \{HC_{n,k}\}, \qquad \mbox{where $HC_{n,k} = \frac{\lambda_k - E_0[\lambda_k]}{SD_0(\lambda_k)}$}. 
\end{equation} 
The current HC is very different from existing versions of HC 
(e.g.,  \cite{DJ04, DJ08}). 

The HC test is an extension of the well-known   Tracy-Widom (TW) test: 
\begin{equation} \label{DefineTW} 
TW_n \equiv   HC_{n, 1}  =  \frac{\lambda_1 - E_0[\lambda_1]}{SD_0(\lambda_1)}.  
\end{equation} 
The TW test is motivated by recent works on the Tracy-Widom law \cite{Tracy-Widom} and has been studied in several recent papers \cite{Johnstone2001, Johnstone2015, Hallin2013, Hallin2014}.   
HC has advantages over TW:  it is known that the empirical eigenvalues are much noisier at the edge than in the bulk, so it is  possible that $HC_{n,1}$ is less powerful than $HC_{n,k}$ for some $k > 1$.  
Numerical results confirm this point; see Section~\ref{sec:Simul}.  

We compare all four tests in Table \ref{tab:4tests}.  Note  that these tests are not tied to the Gaussian model \eqref{Gaussianmodel} and can be implemented in much broader settings. 

\begin{table}[htb]
\caption{Comparison of four test statistics. Both $CS_n^*$ and $HC_n^*$ have a variant,  $CS_n^+$ and $HC_n^+$,  
proposed for the convenience of theoretical study.} \label{tab:4tests}
\scalebox{1.13}{
\begin{tabular}{|c|c|c|c|}
\hline
Trace ($S_n$)    & CuSum ($CS_n^*$) & HC ($HC_n^*$)  & TW ($TW_n$)    \\
\hline
$\sum_{i= 1}^n \lambda_i$    & $\max_{1\leq k \leq n}  \{ \frac{S_k - E_0[S_k]}{SD_0(S_k)}\}$ & $\max_{1\leq k \leq n}  \{ \frac{\lambda_k - E_0[\lambda_k]}{SD_0(\lambda_k)}\}$ & $\frac{\lambda_1 - E_0[\lambda_1]}{SD_0(\lambda_1)}$    \\ [1.05 ex]
\hline
\end{tabular}}
\end{table}

The CS and HC tests are convenient to use in practice, but are {\it difficult} to analyze {\it theoretically}: it is hard to pin down the magnitudes of $SD_0(S_k)$ and $SD_0(\lambda_k)$ analytically. 
Fortunately, we have the following upper bounds \cite{Bao, PillaiYin}. With probability at least $1-o(n^{-1})$, for all $1\leq k\leq n$,  
\begin{equation} \label{regidity} 
SD_0(S_k) \leq (L_n)^{c} (k/n)^{2/3}, \;\;\;\;  
SD_0(\lambda_k) \leq (L_n)^{c} n^{-\frac{2}{3}}[k\wedge(n+1-k)]^{-\frac{1}{3}},  
\end{equation} 
where $L_n\equiv [\log(n)]^{\log(\log(n))}$ and $c > 0$ is a universal constant. These motivate the following {\it variants} of $CS_n^*$ and $HC_n^*$: 
\begin{equation} \label{variants} 
CS_n^+  = \max_{1 \leq k \leq n}\Big\{ \frac{S_k - E_0[S_k]}{(k/n)^{2/3}}\Big\},  \quad    HC_n^+ = \max_{1 \leq k \leq n} 
\Big\{\frac{\lambda_k - E_0[\lambda_k]}{n^{-2/3}[k\wedge(n+1-k)]^{-1/3}}\Big\}.  
\end{equation} 
Due to the explicit forms in the denominators, these variants are more convenient for theoretical analysis than the original version of the statistics.

\subsection{Rare and Weak Spike (RWS) model} 
\label{subsec:RWS} 
Let $X_i \stackrel{iid}{\sim} N(0, \Sigma)$ be as in \eqref{Gaussianmodel}. Our interest is to test whether the null 
\begin{equation} \label{model-null} 
H_0^{(n)}: \qquad \Sigma = I_p  
\end{equation}  
holds or not.  
Fixing an integer $1 \leq r \leq p$ and a parameter $\delta > 0$ and letting $\mathbb{S}(p,r)$ be the Stiefel manifold  \cite{Chikuse} (consisting all matrices $Q \in \mathbb{R}^{p,r}$ such that $Q' Q = I_r$),   
we consider a specific alternative hypothesis: 
\begin{equation} \label{model-alt} 
H_1^{(n)}: \quad  \Sigma = (I_p - \delta QQ')^{-1} = I_p  +  \frac{\delta}{1-\delta} \cdot   Q Q',    
\end{equation}  
where $Q$  is uniformly generated from the Stiefel manifold $\mathbb{S}(p, r)$.

We use $n$ as the driving asymptotic parameter, and tie $(p, r, \delta)$ to $n$ by fixed parameters. In detail,  fixing   $\alpha, \beta   \in (0,1)$ and $\gamma \geq 1$, we assume  
\begin{equation} \label{model-cal2} 
p_n /  n \goto \gamma, \qquad 
r = r_n  = n^{1 - \beta}, \qquad \delta =  \delta_n = n^{-\alpha}.  
\end{equation} 

\bed  
We call (\ref{model-alt})-(\ref{model-cal2}) the Rare and Weak Spike (RWS) model. 
\eed 
Similar Rare/Weak models have been used in many recent works but for different problems \cite{DJ15}, and  the Rare/Weak Spike model  here is new.

{\bf Remark}. ({\it A Gaussian proxy model to RWS}). 
The RWS model is hard to analyze: in the likelihood ratio associated with the testing problem, it is hard to integrate $Q$ out with the law of $Q \sim \mathbb{S}(p,r)$. In the case where $r$ is fixed while $n \goto \infty$, \cite{Hallin2013, Hallin2014} attacked the problem with the so-called Laplace's method and careful large deviation analysis of the spherical integrals, but how to extend their techniques and results to our case (where $r$ grows to $\infty$ at an algebraic rate as $n \goto \infty$)  remains unclear. 

We propose a new technical device by introducing a proxy model that is very close to RWS but is easier to analyze. In detail, we consider a proxy testing problem where we replace 
the alternative hypothesis $H_1^{(n)}$ by   
\[  
\widetilde{H}_1^{(n)}:   \qquad  \widetilde{\Sigma} =  [I_p - (\delta/p)  Y Y']^{-1},  \;\;\;\;  Y = Z \cdot 1\{\|Z\| \leq \frac{1}{2}\sqrt{p/\delta}\}, 
\] 
where $Z \in \mathbb{R}^{p,r}$ is the matrix that has $iid$ $N(0,1)$ entries. Using Le Cam's ``comparison of experiments", we are able to prove that two models are close to each other, for a wide region in the parameter space (see Figure~\ref{fig:phase} and Lemma~\ref{lem:differ}). Note that Gaussian proxy model is comparably easier to analyze than the RWS: the analysis of the former relies on the properties of Wishart matrices where many results exist, while that of the latter relies on
properties of spherical integrals, a topic that is comparably less studied.

\subsection{Main results} 
\label{subsec:main}  
Our result has two parts: phase transition for the testing problem, and bounds for empirical eigenvalues under $H_1^{(n)}$.

In our results on phase transition, especially Theorem \ref{thm:phase}, recent development in Random Matrix Theory (RMT) has played an important role. However, {\it knowledge flow is not a  one-way street}:  our understanding of the testing problem in turn sheds lights on some of the problems in RMT.   

In detail, under  $H_0^{(n)}$, recent works in RMT have shed interesting   lights on the bounds of the empirical eigenvalues (e.g., \cite{Yau}), but under   $H_1^{(n)}$,  such bounds  are  much less studied and remain  largely unknown.  
Interestingly, Theorem \ref{thm:phase} provides an approach to studying the bounds under $H_1^{(n)}$.   

The idea is that, for parameters in the {\it Region of Impossibility} (see below),  any test is asymptotically powerless in distinguishing $H_1^{(n)}$ from $H_0^{(n)}$. We can therefore 
use existing bounds on the empirical eigenvalues under $H_0^{(n)}$ to derive similar bounds  under $H_1^{(n)}$; the resultant bounds are non-trivial to derive using RMT, at least for some of the parameter ranges.  


\begin{figure}[tb!]  \label{fig:phase} 
\centering
\includegraphics[width=0.45\textwidth]{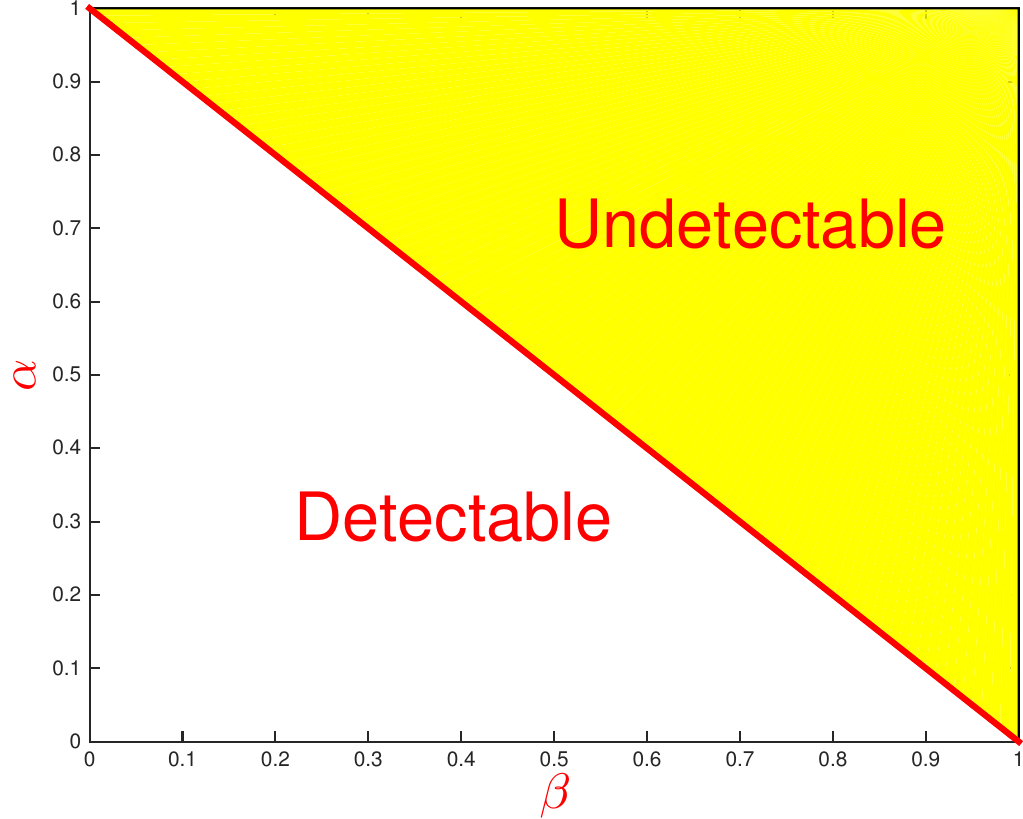} $\;$
\includegraphics[width=0.44\textwidth]{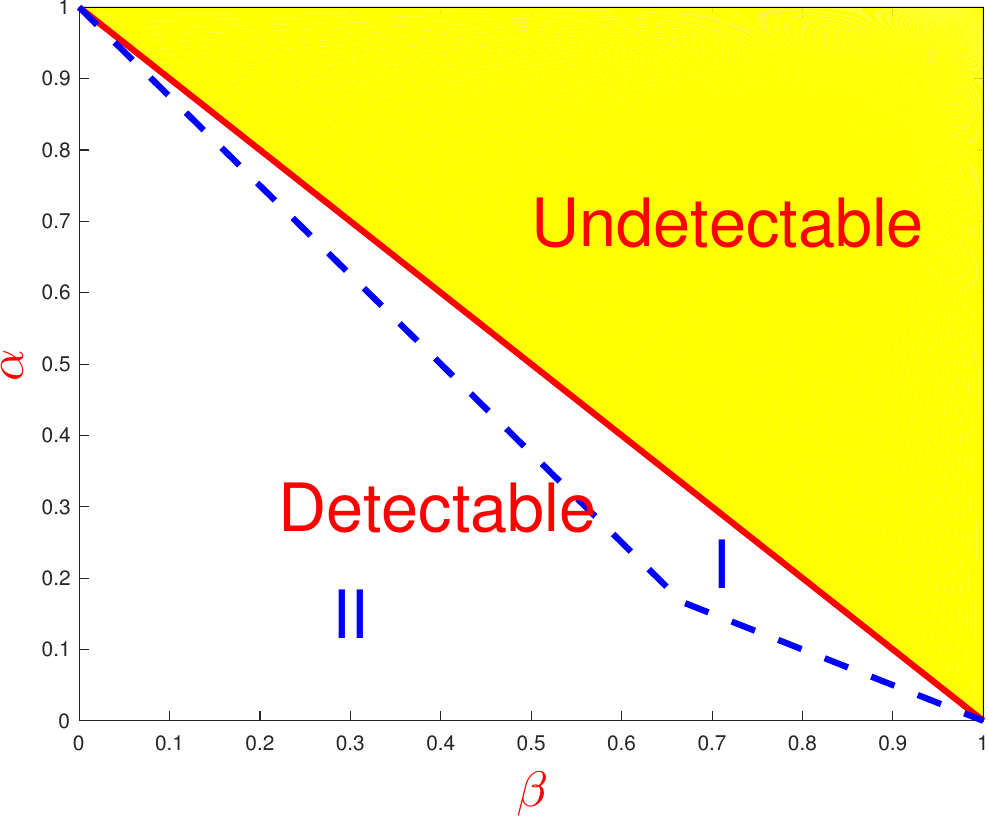}
\caption{Left:  the curve $\alpha = (1 - \beta)$ partitions the  phase space into the ``detectable region" and the ``undetectable region".  Right: the curve $\alpha=\max\{1-5\beta/4, (1-\beta)/2\}$   (see Lemma~\ref{lem:differ}) partitions the ``detectable region" into part I and II.   In both the undetectable region and part I of the detectable region,  the RWS model and the Gaussian proxy model have negligible differences.}  \label{fig:phase}
\end{figure} 

Now, first, consider the trace-based test statistic $S_n$. Under $H_0^{(n)}$,  
$nS_n \sim  \chi_{np}^2(0)$, and so  $(2np)^{-1/2} (nS_n - np) \approx N(0,1)$ for large $(n, p)$.  Fixing  a parameter $q > 0$,  suppose we reject $H_0^{(n)}$  if and only if 
\begin{equation} \label{tracetest} 
(2np)^{-1/2} (nS_n - np) \geq \sqrt{2 q \log(n)}, \qquad  q > 0. 
\end{equation} 
The following theorem establishes the phase transition associated with the testing problem (\ref{model-null})-(\ref{model-alt}) and 
the optimality of the trace-based test $S_n$. 
\begin{thm} \label{thm:phase}  
({\it Phase transition and optimality of the trace-based test}). 
Fix $\alpha, \beta \in (0,1)$,  $\gamma \geq 1$, and $q > 0$.  Suppose $X_i \stackrel{iid}{\sim} N(0, \Sigma)$  as in \eqref{Gaussianmodel}. Consider the testing problem \eqref{model-null}-\eqref{model-alt}
 where $\Sigma$ satisfies the RWS model \eqref{model-alt}-\eqref{model-cal2}. 
\begin{itemize} 
\item If $\alpha + \beta > 1$, then as $n \goto \infty$,  for any test, the sum of Type I and Type II testing errors   tends to $1$. 
\item If $\alpha + \beta < 1$, then as $n \goto \infty$,  the level of the trace-based test in \eqref{tracetest} tends to $0$ while the power tends to $1$.
\end{itemize}   
\end{thm} 
See Figure \ref{fig:phase}.  The proof of the second bullet point is not hard so is omitted.  The first bullet point is proved in Section \ref{sec:main},  where the key ingredients include (a) justifying that the RWS is close to the Gaussian proxy model, and (b) analyzing the  proxy model with delicate Random Matrix Theory. 
Note that the case $\alpha + \beta = 1$ is more delicate and is addressed in Section \ref{subsec:criticalbdry}.

Second, we consider the CS test statistic, focusing on $CS_n^+$ (the variant of $CS_n^*$).  
By definition of $CS_n^+$ (see \eqref{regidity}-\eqref{variants}), we reject $H_0^{(n)}$ if and only if 
\begin{equation} \label{CScritical} 
CS_n^+ \geq \tilde{L}_n,  \qquad \mbox{where } \tilde{L}_n=  [\log(n)]^{\log^2(\log(n))}. 
\end{equation} 
Recall that $S_k = \lambda_1 + \ldots + \lambda_k$ is the $k$-th cumulative sum of the empirical eigenvalues.   
The following theorem shows that  $CS_n^+$ is an optimal test and  provides some non-trivial bounds on $S_k$ under the alternative. 
\begin{thm} \label{thm:CS}  
({\it Optimality of the CuSum test and bounds for $S_k$}). 
Fix $\alpha, \beta \in (0,1)$ and  $\gamma \geq 1$.  Suppose $X_i \stackrel{iid}{\sim} N(0, \Sigma)$  as in \eqref{Gaussianmodel} and let $\tilde{L}_n=[\log(n)]^{\log^2(\log(n))}$. Consider the testing problem \eqref{model-null}-\eqref{model-alt} where $\Sigma$ satisfies the RWS model \eqref{model-alt}-\eqref{model-cal2}. 
\begin{itemize} 
\item If $\alpha + \beta  <  1$, then as $n \goto \infty$, the level of the CS test defined in \eqref{variants} and \eqref{CScritical} tends to $0$ while the power tends to $1$.  
\item If $\alpha + \beta > 1$ and $H_1^{(n)}$ holds, then for sufficiently large $n$, 
with probability at least $1 - C\log(n)n^{1 - (\alpha + \beta)}$, 
\begin{equation} \label{Skbd} 
\bigl|S_k - E_0[S_k] \bigr| \leq C \tilde{L}_n (k/n)^{2/3},  \qquad \mbox{for all $1 \leq k \leq n$}.
\end{equation} 
\end{itemize} 
\end{thm}  
To the best of my knowledge,  the bounds in (\ref{Skbd}) are new.  
The probability bound on the excluded event (i.e.,  $C\log(n)n^{1 - (\alpha + \beta)}$)  is derived from the relationship between the $L^1$-distance and the Neyman-Pearson Lemma \cite{Aad}; the bound may be improved, using presumably a different technique.

Next,  consider the HC test statistic, focusing on $HC_n^+$ (the variant of $HC_n^*$).   
By definitions (i.e., \eqref{regidity}-\eqref{variants}), we reject $H_0^{(n)}$ if and only if 
\begin{equation} \label{HCcritical} 
HC_n^+ \geq \tilde{L}_n, \qquad \mbox{where } \tilde{L}_n=  [\log(n)]^{\log^2(\log(n))}. 
\end{equation} 
\begin{thm} \label{thm:HC}  
({\it Behavior of the HC test and bounds for $\lambda_k$}). 
Fix $\alpha, \beta \in (0,1)$ and $\gamma >  1$.\footnote{The case of $\gamma = 1$ is more complicated, for the behavior of the smallest eigenvalues is different. To save space, we omit discussions of this case.}   Suppose $X_i \stackrel{iid}{\sim} N(0, \Sigma)$  as in \eqref{Gaussianmodel}.  Let $\tilde{L}_n=[\log(n)]^{\log^2(\log(n))}$. Consider the testing problem \eqref{model-null}-\eqref{model-alt} where $\Sigma$ satisfies the RWS model \eqref{model-alt}-\eqref{model-cal2}. 
\begin{itemize} 
\item If $\alpha + \beta  <  1$, then  as $n \goto \infty$, the level of  the HC  test defined in \eqref{variants} and \eqref{HCcritical} tends to $0$. 
\item If $\alpha + \beta   > 1$ and $H_1^{(n)}$ holds, then for sufficiently large $n$, with probability at least $1 - \log(n) n^{1 - (\alpha + \beta)}$, 
\[
\bigl|\lambda_k - E_0[\lambda_k] \bigr| \leq \tilde{L}_n n^{-2/3} [k\wedge(n+1-k)]^{-1/3}, \qquad \mbox{for all $1 \leq k \leq n$}.  
\] 
\end{itemize} 
\end{thm} 
Compared to bounds of $\lambda_k$ in literature (e.g., \cite{BAP, Baik, Paul2007}),  our results are new for (a) the literature focus on the case where $r_n$ is fixed as $n \goto \infty$, while our results are for  the case of $r_n = n^{1 - \beta}$, and (b) the literature focus on bounds for only edge eigenvalues   and our bounds are for all eigenvalues.     

Last, consider the Tracy-Widom test statistic  $TW_n$. It is known that under $H_0^{(n)}$, $TW_n$ converges weakly to the Tracy-Widom law. Let $F_{TW}$ be the CDF of the Tracy-Widom law.  It is known that $1-F_{TW}(t)\sim e^{-(2/3)t^{3/2}}$ \cite{dumaz2013} for large $t$. In light of this, we reject $H_0^{(n)}$ if and only if 
\begin{equation} \label{TWcritical} 
TW_n \geq [3\log(n)]^{2/3}.  
\end{equation} 
 
\begin{cor} \label{cor:TW}  
({\it Sub-optimality of the TW test}). 
Fix $\alpha, \beta \in (0,1)$ and  $\gamma  \geq   1$. Suppose $X_i \stackrel{iid}{\sim} N(0, \Sigma)$ as in \eqref{Gaussianmodel}. Consider the testing problem \eqref{model-null}-\eqref{model-alt} where $\Sigma$ satisfies the RWS model \eqref{model-alt}-\eqref{model-cal2}.  If $\alpha + \beta  <  1$, then  as $n \goto \infty$, the level of the TW  test defined in \eqref{HCcritical} tends to $0$. If additionally $\alpha > 2/3$, then the power of the TW test also tends to $0$.    
\end{cor}  
The TW test is asymptotically powerless 
in the sub-region 
\[
\{(\beta, \alpha):  0 < \beta < 1,  \alpha > 2/3, \alpha + \beta < 1\}; 
\]
In this region, both the CS test and the trace-based test have asymptotically full power.  
This result is consistent with our numerical study (see Section~\ref{sec:Simul}), where we observe that 
when there are multiple spikes, the TW test usually behaves unsatisfactorily.  
 
 
{\bf Remark}. ({\it Comparison with the testing problem of sparse normal means}).  
Consider a setting where we have samples $y_i \sim N(\mu_i, 1)$, $1 \leq i \leq n$. We are 
interested in testing whether $\mu_i = 0$ for all $1 \leq i \leq n$, or that most of $\mu_i$ are  $0$ except for a small fraction of them are nonzero. 
In such a sparse normal means test setting, it is preferable to use a small fraction of extreme observations (instead of the bulk) for testing (e.g., \cite{DJ04}). 

\begin{figure}[tb!]  
\centering
\includegraphics[width=0.4\textwidth, height=.25\textwidth]{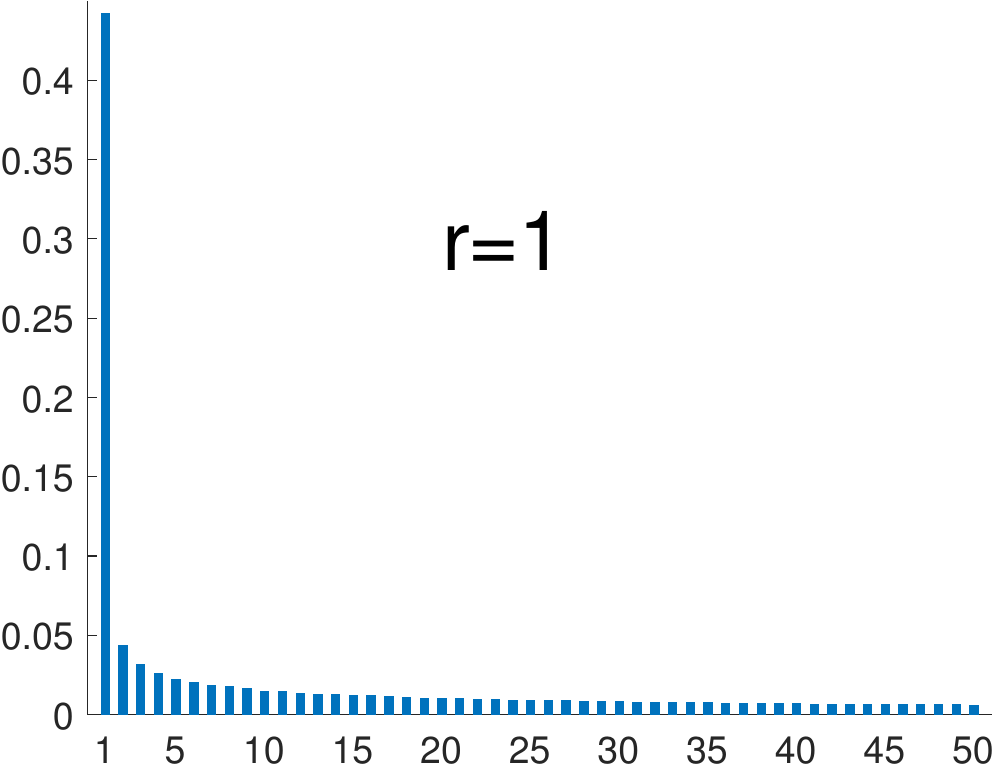}
\includegraphics[width=0.4\textwidth, height=.25\textwidth]{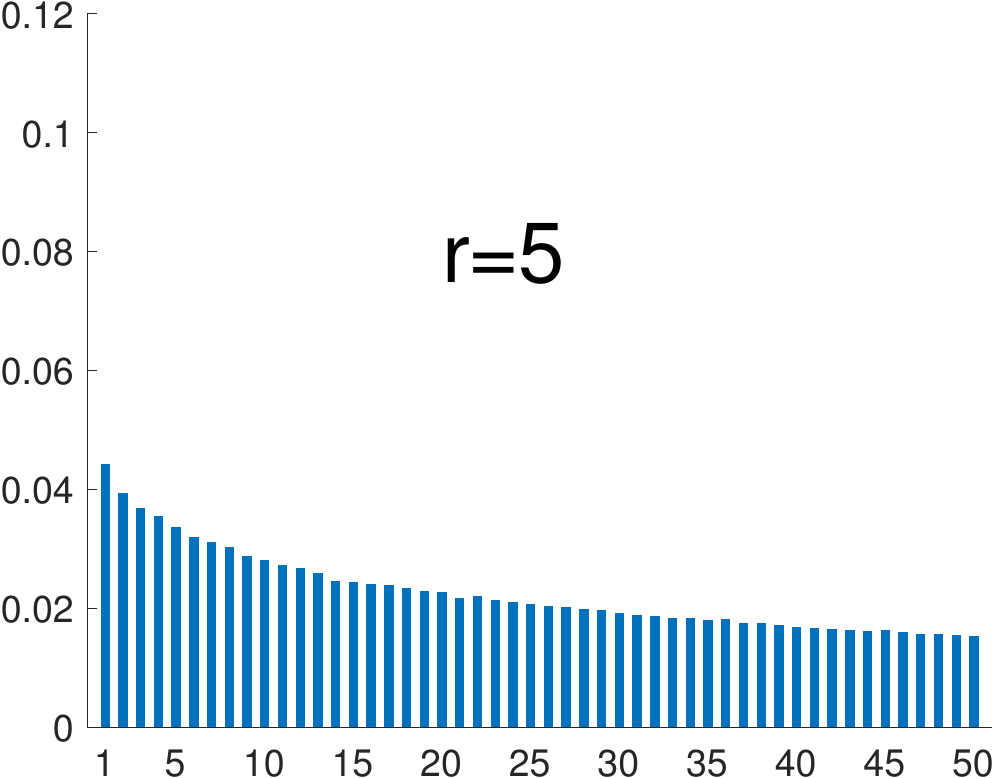}
\caption{Plots of $E_1[\lambda_k]-E_0[\lambda_k]$ for $k=1,\cdots, 50$, where $(n,p)=(1000, 1200)$, $(r,\delta)=(1,2)$ in left panel, and $(r,\delta)=(5, 0.8)$ in right panel. } \label{fig:sparse}
\end{figure}

To those familiar with the normal means problem, Theorem \ref{thm:phase}  may strike as a surprise, for it shows the trace-based test is optimal, even when the spikes are very sparse; this is so because the current setting is very different from the normal means settings:   
every spike  in $\Sigma$ affects the {\it bulk} of $\lambda_i$ in a subtle way; and even when the spikes in $\Sigma$ are very sparse, the sparsity is lost when we 
look at the vector 
\begin{equation} \label{vector_add} 
\big(E_1[\lambda_1] - E_0[\lambda_1], \; E_1[\lambda_2] - E_0[\lambda_2],\; \ldots,\; E_1[\lambda_n] - E_0[\lambda_n] \big),  
\end{equation} 
where $E_0$ and $E_1$ are the expectation under $H_0^{(n)}$ and $H_1^{(n)}$, respectively. In Figure~\ref{fig:sparse}, we plot the vector in (\ref{vector_add})  for the cases $r=1$ and $r=5$ ($n,p,\delta$ are the same as those of Table~\ref{tab:HC}). We observe that as soon as $r$ move away from $1$, the vector in (\ref{vector_add}) become resonably non-sparse. 
 
{\bf Remark}. ({\it Other trace-based statistics}). 
One possible variant of the trace-based test $S_n$ is $S_n^{(2)}=\sum_{k=1}^n (\lambda_k-1)^2$. 
The test targets on the case where the off-diagonals of $\hat{\Sigma}$ contain stronger evidence  against the null  than the diagonals,  so it is not surprising that,  asymptotically,  for $(\alpha, \beta)$ such that  $2 \alpha + \beta > 1$ while $\alpha + \beta < 1$,  $S_n^{(2)}$ is  powerless while $S_n$ has full power.  On the other hand, $S_n^{(2)}$ may be more powerful in some other cases (e.g., \cite{ChenSX}).

%

\subsection{Extensions and Le Cam's comparison of experiments} 
\label{subsec:exten} 
In his work on comparison of experiments, Le Cam asserts that ``adding noise always makes the inference more difficult".  This allows us to compare our setting with many other settings, and generalize our lower bound argument in Theorem \ref{thm:phase} (one of our major contribution in this paper) to much broader settings. 


We now consider two extensions of our lower bound argument. 
In the first extension,  for two covariance matrices $\Sigma$ and $\widetilde{\Sigma}$, 
we compare two experiments.  In the first one, we are interested in testing 
\[
H_0^{(n)}: X_i  \stackrel{iid}{\sim} N(0, \Sigma)   \qquad vs.   \qquad H_1^{(n)}: (X_i |Q) \stackrel{iid}{\sim} N(0, \Sigma + \frac{\delta}{1 - \delta}  QQ'),  
\]
and in the second one, we are interested in testing 
\[
\widetilde{H}_0^{(n)}: X_i  \stackrel{iid}{\sim} N(0, \widetilde{\Sigma})   \qquad vs.   \qquad \widetilde{H}_1^{(n)}: (X_i |Q) \stackrel{iid}{\sim} N(0,  \widetilde{\Sigma} + \frac{\delta}{1 - \delta}  QQ'), 
\]
where as before, $Q$ is uniformly distributed over the Stiefel manifold $\mathbb{S}(p,r)$. 
Let $f_0, f_1, \tilde{f}_0, \tilde{f}_1$ be the joint density of $\{X_i\}_{i = 1}^n$ under $H_0^{(n)}$, $H_1^{(n)}$, $\widetilde{H}_0^{(n)}$, and $\widetilde{H}_1^{(n)}$, respectively. We assume $\Sigma  \preceq \widetilde{\Sigma}$, so the second experiment can be viewed as the result of adding noise to the first experiment. Lemma \ref{lemma:LeCam2} solidifies the claim that  ``adding noise makes the testing problem harder".  
\begin{lemma} \label{lemma:LeCam2} 
If $\Sigma \preceq \widetilde{\Sigma}$, then $\|\tilde{f}_1 - \tilde{f}_0\|_1 \leq \|f_1 - f_0\|_1$.  
\end{lemma}
Applying Lemma~\ref{lemma:LeCam2} with $\Sigma = c_0 I_p$  and $\widetilde{\Sigma} =  \Sigma^*$ (see below), we have the following theorem, as a direct result of 
Theorem~\ref{thm:phase}. 
\begin{thm} \label{thm:LeCam2} 
({\it Extension of lower bound argument, I}). 
Fix $\alpha, \beta\in (0,1)$, $c_0 > 0$,  and   $\gamma \geq 1$.   
Given $n$ independent vectors $X_i \in \mathbb{R}^p$, $1 \leq i \leq n$,  we are interested in testing $H_0^{(n)}: X_i \stackrel{iid}{\sim} N(0, \Sigma^*)$ versus the alternative $H_1^{(n)}:  (X_i | Q) \stackrel{iid}{\sim}  N(0,  \Sigma^* + \frac{\delta}{1 - \delta}  QQ')$, where $Q$ is uniformly distributed 
over the Stiefel manifold   $\mathbb{S}(p,r)$, $\Sigma^* \in \mathbb{R}^{p,p}$ is a covariance matrix, and  $(\delta, r, p)$ satisfies (\ref{model-cal2}) as in Theorem \ref{thm:phase}. If 
$\alpha + \beta > 1$ and $\lambda_{\min}(\Sigma^*)\geq c_0$,\footnote{For a symmetric matrix $A$, $\lambda_{\min}(A)$ stands for the minimum eigenvalue of $A$.} then as $n \goto \infty$, the sum of Type I and Type II errors of any test tends to $1$.  
\end{thm} 

For another extension of the lower bound argument,   
we show that ``either adding more spikes or increasing the spike strengths makes the testing problem   (\ref{model-null})-(\ref{model-alt})  easier".  Consider two experiments. In the first one, we are interested in testing 
\[
H_0^{(n)}: X_i  \stackrel{iid}{\sim} N(0, I_p)   \qquad vs.   \qquad H_1^{(n)}: (X_i |Q) \stackrel{iid}{\sim} N(0, I_p + Q D Q').
\]
In the second one, we are interested in testing 
\[
H_0^{(n)}: X_i  \stackrel{iid}{\sim} N(0, I_p)   \qquad vs.   \qquad \widetilde{H}_1^{(n)}: (X_i |Q) \stackrel{iid}{\sim} N(0,  I_p + Q \widetilde{D} Q'), 
\]
where $D, \widetilde{D} \in \mathbb{R}^{r,r}$ are positive semi-definite matrices,  and  $Q$ is uniformly distributed over the Stiefel manifold $\mathbb{S}(p,r)$.  
Let $f_0, f_1,   \tilde{f}_1$ be the joint density of $\{X_i\}_{i = 1}^n$ under $H_0^{(n)}$, $H_1^{(n)}$,   and $\widetilde{H}_1^{(n)}$, respectively.  Similarly, we have the following lemma, which is proved  in Section~\ref{sec:main} with nontrivial efforts. 
\begin{lemma} \label{lemma:LeCam3} 
If $D \preceq \widetilde{D}$, then $\|f_1 - f_0\|_1\leq \|\tilde{f}_1 - f_0\|_1$.  
\end{lemma} 
We also have the following theorem. 
\begin{thm} \label{thm:LeCam3} 
({\it Extension of lower bound argument, II}). 
Fix $\alpha, \beta \in (0,1)$, $c_0 > 0$,  and  $\gamma \geq 1$.   
Given $n$ independent vectors $X_i \in \mathbb{R}^p$, $1 \leq i \leq n$,  we are interested in testing $H_0^{(n)}: X_i \stackrel{iid}{\sim} N(0, I_p)$ versus the alterative $H_1^{(n)}:  (X_i | Q) \stackrel{iid}{\sim}  N(0, I_p +  QD^*Q')$, where $Q$ is uniformly distributed 
over the Stiefel manifold   $\mathbb{S}(p,r)$, $D^* \in \mathbb{R}^{r,r}$ is a positive semi-definite matrix, and  $(\delta, r, p)$ satisfies (\ref{model-cal2}) as in Theorem \ref{thm:phase}. If 
$\alpha + \beta > 1$ and $\|D^*\| \leq c_0 \frac{\delta}{1-\delta}$, then as $n \goto \infty$, the sum of Type I and Type II errors of any test tends to $1$.   
\end{thm} 
For a proof,  assume $c_0 = 1$ without loss of generality. 
Applying Lemma \ref{lemma:LeCam3} with $D = D^*$  and $\widetilde{D} =  \frac{\delta}{1 - \delta} I_r$,  
Theorem \ref{thm:LeCam3} follows from Theorem \ref{thm:phase}.

\subsection{Spike detection in the critical case of $r_n \delta_n = O(1)$} 
\label{subsec:criticalbdry}  
In Theorem \ref{thm:phase}, we have either $r_n  \delta_n \goto \infty$ or  $r_n \delta_n \goto 0$. 
We consider a more subtle case where $r_n \delta_n = O(1)$ (in RWS, this is the case of 
$\alpha + \beta =1$). 

For analysis,  we continue to use RWS, except for  
that  $(r_n, \delta_n)$ are calibrated slightly differently.  
In detail, fixing the parameter $\beta\in(0,1)$, $\gamma\geq 1$ and $\theta > 0$, we assume as $n \goto \infty$, 
\begin{equation} \label{newcalibration} 
p_n/n\goto \gamma, \qquad r=r_n=n^{1-\beta},  \qquad  \delta= \delta_n =r_n^{-1} (\theta \sqrt{2 \gamma}). 
\end{equation} 
The following Theorem is proved in Section \ref{sec:main}.  
\begin{thm} \label{thm:critical1} 
({\it Critical case of $\alpha+\beta=1$}). 
Fix $\beta\in (0,1)$, $\theta>0$, and $\gamma\geq 1$. 
Consider the testing problem \eqref{model-null}-\eqref{model-alt}
where \eqref{newcalibration} holds.
As $n\goto \infty$, the log-likelihood ratio $\log(LR_n)$ converges weakly to $N(\mp  \frac{\theta^2}{2},  \theta^2)$, under the null and under the alternative, respectively.  
\end{thm} 
Recall that $S_n$ denotes the trace-based test statistic. By elementary statistics, it is seen that (the convergence is weak convergence): 
\[
\theta\frac{(nS_n - np)}{\sqrt{2np}} -\frac{\theta^2}{2}\goto  N(\mp \frac{\theta^2}{2}, \theta^2),  \;\;\;   \mbox{under $H_0^{(n)}$ and $H_1^{(n)}$, respectively}.   
\] 
This suggests that the trace-based statistic $S_n$ is asymptotically efficient. 
 
In our case, $r_n \delta_n = O(1)$ but $r_n \goto \infty$. A closely related case is that 
as $n \goto \infty$, both $(r, \delta)$ are fixed. Such a case  was studied in detail in \cite{Hallin2014}, with very different techniques. 
Our framework can be extended to such a case. Let $\psi(\lambda)=\frac{\gamma}{\delta}\log(1+\frac{\delta}{\gamma(1-\delta)} - \frac{\delta}{\gamma}\lambda)$. Let $\mu_{\hat{\Sigma}}$ be the empirical spectral measure associated with $\hat{\Sigma}$ and let $\mu_{n,p}$ be the Marchenko-Pastur law (see \eqref{MP} for definition). 
By a simple modification of the proof of Theorem~\ref{thm:critical1}, we can show that
$\log(LR_n)=-\frac{nr\delta}{2}[\int g(\lambda)\mu_{\hat{\Sigma}}(d\lambda) -\int g(\lambda)\mu_{n,p}(d\lambda) ] - \frac{r^2\delta^2}{2\gamma(1-\delta)}+o_P(1)$. 
So $\log(LR_n)$ converges to a weak limit according to the central limit theorem for linear spectral statistics \cite{BaiSilverstein}. It yields that $\log(LR_n)\goto N(\mp\tfrac{\tilde{\theta}^2}{2}, \tilde{\theta}^2)$ under $H_0^{(n)}$ and $H_1^{(n)}$, respectively,  where $\tilde{\theta}=r[-\tfrac{1}{2}\log(1-\tfrac{\delta^2}{\gamma})]^{1/2}$. This result coincides with that of \cite[Proposition 4]{Hallin2014}.

\subsection{A stylized application: Covert Communications} 
\label{subsec:covert} 
We wish to communicate with our ``friends" through a covert channel, and our ``enemy" is trying to intercept it.  We encode the desired information by a length-$p$ string  of  ``0" and ``1", denoted by 
$\eta$ (note $\eta$ has $r$ nonzero entries).  Fix $\delta > 0$ and a $p \times p$ matrix $\widetilde{Q} = [q_1, q_2, \ldots, q_p]$.   
We send our friends vectors $X_i \in \mathbb{R}^p$   
\[
X_i =  \delta  \sum_{\{k: \eta(k) = 1 \}} z_k  \cdot   q_k  + Z_i, \qquad 1 \leq i \leq n, 
\]
where $z_k \stackrel{iid}{\sim} N(0,1)$,  $Z_i \stackrel{iid}{\sim} N(0, I_p)$, and they are independent.    
Let $Q$ be the $p \times r$ matrix $\{q_k: \eta(k) = 1, 1 \leq k \leq p\}$. Note that $X_i$ can be equivalently viewed as samples from 
$N(0, \Sigma)$,  where $\Sigma = I_p + \delta QQ'$.  

Assume that the matrix $\widetilde{Q}$ is available to our ``friends", but not to our ``enemy". 
When our ``friends" receive $X_i$,  first,  they can obtain $\widetilde{X}_i$:  
\[
\widetilde{X}_i = \delta \sum_{\{k: \eta(k) = 1\}} z_k  \cdot  e_k   + \widetilde{Z}_i,    \qquad \widetilde{X}_i = \widetilde{Q}' X_i,  \qquad  \widetilde{Z}_i = \widetilde{Q}' Z_i.  
\]
Next, from $\widetilde{X}_i$,   they are able to retrieve the vector $\eta$, provided with some mild conditions on $(r, \delta)$.  
For our ``enemy", $\widetilde{Q}$ is not available, so the problem reduces to our previous setting of (\ref{model-alt}).   Note also that if $\widetilde{Q}$ is randomly generated, then $Q$ is uniformly distributed in the Stiefel manifold $\mathbb{S}(p, r)$. Applying Theorem~\ref{thm:phase} and conventional results of sparse normal-means problems, we have the following theorem:
\begin{thm} \label{thm:covert}
({\it Covert communitcation}).
Fix $\alpha, \beta \in (0,1)$ and $\gamma \geq 1$.  Suppose  $(\delta, r,  Q)$ satisfy model (\ref{model-alt})-(\ref{model-cal2}), and where  $1/2 < \beta < 1$ and $(1 - \beta) < \alpha < 1/2$. Then with high probability,   our ``friends" are able to exactly decode the vector $\eta$,  while our ``enemy" is not able to even distinguish whether we are sending some signals or we are merely sending white noise. 
\end{thm} 


\subsection{Summary} 
\label{subsec:summary} 
Our results provide both a better understanding of the problem of detecting rare/weak spikes, and 
better eigenvalue bounds  for the empirical covariance matrix $\hat{\Sigma}$ (based on $iid$ samples from $N(0, \Sigma)$).  
 
For the testing problem, we adopt a Rare and Weak Spike (RWS) model and discover an interesting phase transition. We study four different tests, where the Higher Criticism test and the CuSum test are new. We show that both the trace test and the CuSum test are optimal. We find that even when the spikes are very rare, it is still necessary to use {\it many} eigenvalues (instead of a few extreme eigenvalues) 
for testing against the null; such a finding is very different from those in testing settings regarding sparse normal means. 

Motivated by the interest of ``$p>n$" in modern applications, we assume $\gamma\geq 1$, but this should not be taken as a constraint. When $p/n\to\gamma<1$, the detection boundary is the same, and the asymptotic behavior of all four tests can be studied similarly. 

Our study is connected to the large body of literature on testing of sphericity \cite{bai2009corrections, BR13, cai2013optimal, ChenSX, dobriban2016sharp, Johnstone2001, Johnstone2015, ledoit2002some, Hallin2013, Hallin2014}, but most of these works focus on specific tests,  not on the phase transition.  
A few exceptions are \cite{BR13, cai2013optimal, Hallin2013, Hallin2014}, but the alternative hypotheses   considered there are very different from ours: in our setting, the number of spikes grow rapidly with $n$;   such a case has not been studied in the literature. 

In \cite{Hallin2013,Hallin2014}, they considered a testing framework equivalent to our RWS model with both $(r,\delta)$ being constants. The test statistic they proposed is a function of all empirical eigenvalues, resulting from a Laplace approximation of the log-likelihood ratio. They showed that the asymptotic power of their test is better than that of the Tracy-Widom test. In a high level, both their works and our work point out the advantage of using bulk eigenvalues in the test, even for a small $r$. However, there are several major differences between the two works: We focus on the setting of Rare/Weak Spikes where $r\to\infty$ and $\delta\to 0$; there is no straightforward extension of their method/theory to our setting. Their test is essentially the likelihood ratio test, which requires knowing $(r,\delta)$, but all four tests considered in our work are adaptive, whose construction doesn't depend on parameters of the alternative. The technical approaches are also different. Their main tools are the Laplace's method and large-deviation analysis of spherical integrals; how to adapt these techniques to our setting is unclear. We use the Gaussian proxy model and Le Cam's comparison of experiments. The Gaussian proxy model and the theory of its intimacy to RWS (Lemma~\ref{lem:differ}), as a new technical device, will be useful for studying other problems associated with RWS.  

For bounds on the eigenvalues of $\hat{\Sigma}$,  while most of the literature have been focused on the case 
of $\Sigma = I_p$, our results focus on the case where $\Sigma$ is a spike matrix.  
For the latter, a few works exist (e.g.,  \cite{BAP,Baik, Paul2007}), but the focus there is on the bounds for the extreme eigenvalues {\it only},  and is for the case of {\it finitely} many spikes. 
In comparison,  our bounds are for all eigenvalues and are for the case where the number of spikes grows rapidly with $n$ as $n \goto \infty$.  See Section~\ref{sec:discuss} for more discussion. 

Our philosophy is that, knowledge flow is a two-way street: while RMT may help us obtain better results in statistical inference, statistical inference can also provide better results in RMT,  in return.   


\subsection{Content and notations} 
\label{subsec:content} 
The remaining part of the paper is organized as follows. 
Section~\ref{sec:Simul} contains a small-scale numerical study. In Section~\ref{sec:RMT}, we present a list of useful results on Random Matrix Theory (RMT), including some new results; this section can be read independently. In Section~\ref{sec:main}, we prove the main theorems and lemmas. Section~\ref{sec:discuss} contains discussions. Proofs of the secondary lemmas are relegated to the appendix. 

For two real numbers $a$ and $b$, $a\wedge b$ and $a\vee b$ denote the minimum and maximum of them, respectively. We say two positive sequences $a_n\sim b_n$ if $a_n/b_n\goto 1$ as $n\goto\infty$, $a_n\asymp b_n$ if $|a_n/b_n|$ are uniformly upper and lower bounded by constants. 
When $\xi$ is a vector, $\|\xi\|$ denotes the vector $L^2$-norm; when $\xi$ is a matrix, $\|\xi\|$ denotes the matrix spectral norm and $\|\xi\|_F$ denotes the matrix Frobenius norm. For two square matrices $A,B$, we say $A\preceq B$ 
 if $B-A$ is positive semidefinite. For two probability densities $f$ and $g$,  $\|f-g\|_1\equiv \int |f(x)-g(x)|dx$ denotes the $L^1$-distance between them.

\section{Numerical study} 
\label{sec:Simul}  

We investigate the four tests with a small-scale experiment. Fix $(n, p) = (1000,1200)$. 
For each $(r, \delta)$, we generate $100$ data sets of $\{X_i\}_{i = 1}^n$ from \eqref{model-null} to represent the null, 
and generate $100$ data sets of $\{X_i\}_{i=1}^n$ from \eqref{model-alt} to represent the alternative, and apply all four tests. 
We measure the ``ideal testing error'' for each test which corresponds to the rejection threshold that minimizes the sum of type I and type II errors over $200$ data sets. For each setting, we run $20$ repetitions. 

First, we compare the HC test with the TW test.\footnote{For HC, $E_0[\lambda_k]$ and $SD_0(\lambda_k)$ are computed by simulating the null for $10,000$ times. Same for CuSum.} The HC statistic is the maximum of $HC_{n,k}$ over $k$, and each $HC_{n,k}$ 
can be used as a test statistic, including TW as a special case with $k = 1$. We record the ``ideal testing error" of $HC_{n,k}$ for all $1\leq k\leq n$, and let $\hat{k}_{hc}$ minimize this error. 
We also record the ``ideal testing error" of the HC test, where $k$ is chosen adaptively by data. The results are summarized in Table~\ref{tab:HC}.
\begin{table}[htb!]
\caption{Comparison of TW and HC tests (standard deviations are in brackets).} \label{tab:HC}
\scalebox{1}{
\begin{tabular}{|c||c||c|c|}
\hline 
$(r,\delta)$& $\hat{k}_{hc}$ & TW   & HC \\
\hline
$(1,2)$  &  $1(0)$ & $0(0)$   & $.001(.003)$  \\ 
$(5, .8)$  & $22 (9.7)$ & $.48(.05)$  & $.28(.04)$     \\ 
$(10, .5)$  &  $53(25)$ & $.56(.05)$  & $.19(.03)$ \\ 
\hline
\end{tabular}}
\end{table}

The results suggest that (a) when $r$ is very small, $\hat{k}_{hc}$ is very close to $1$, and two tests have very similar 
behaviors, with TW being slightly better, and (b) when $r$ gets larger, $\hat{k}_{hc}$ is bounded away from $1$, and the HC test is usually better, and significantly so when $r$ increases.

Next, we compare the trace test, CuSum test, and TW test. Each cumulative sum $CS_{n,k}$ can be used as a test statistic, with TW and trace being two special cases of $k=1$ and $k=n$, respectively. Similarly, we let $\hat{k}_{cs}$ be the $k$ such that the ``ideal testing error" of $CS_{n,k}$ is minimized. The results are summarized in Table~\ref{tab:compare4tests}. 
\begin{table}[htb!]
\caption{Comparison of TW, Trace, and CuSum (standard deviations are in brackets).} \label{tab:compare4tests}
\scalebox{1}{
\begin{tabular}{|c||c||c|c|c|}
\hline 
$(r,\delta)$   &   $\hat{k}_{cs}$ & TW  & Trace   & CuSum       \\
\hline
$(1,2)$  &   $1(0)$ & $0(0)$ & $.48(.05)$ & $0(0)$ \\ 
$(5, .8)$  &  $125(65)$ & $.48(.05)$ &  $.18(.03)$ &  $.16(.03)$    \\ 
$(10, .5)$ & $321(207)$ & $.56(.05)$ & $.08(.03)$  & $.12(.03)$   \\ 
\hline
\end{tabular}}
\end{table}

It suggests that (a) when $r$ is very small,  $\hat{k}_{cs}$ is very close to $1$, and TW and CuSum have similar 
performances, and (b) when $r$ gets larger, $\hat{k}_{cs}$ is much larger than $1$, but it is also smaller than $n$; CuSum performs better than the trace test, and the trace test outperforms the TW test.

\begin{figure}[tb!] 
\centering
\includegraphics[width=0.4\textwidth]{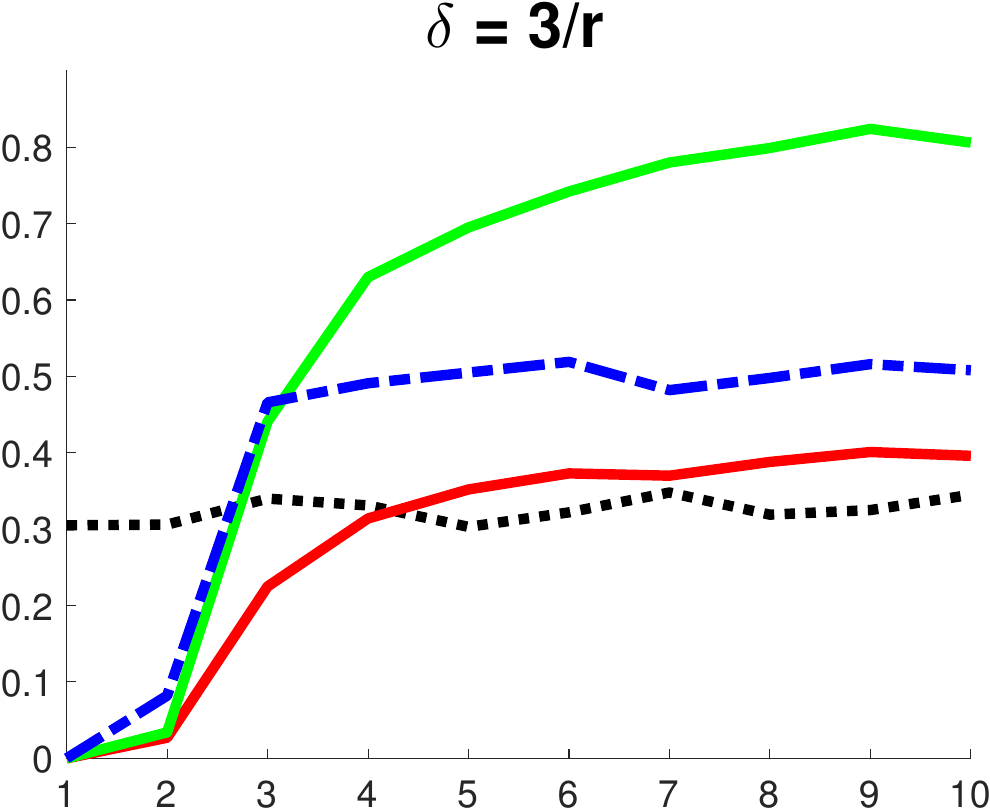}
\includegraphics[width=0.4\textwidth]{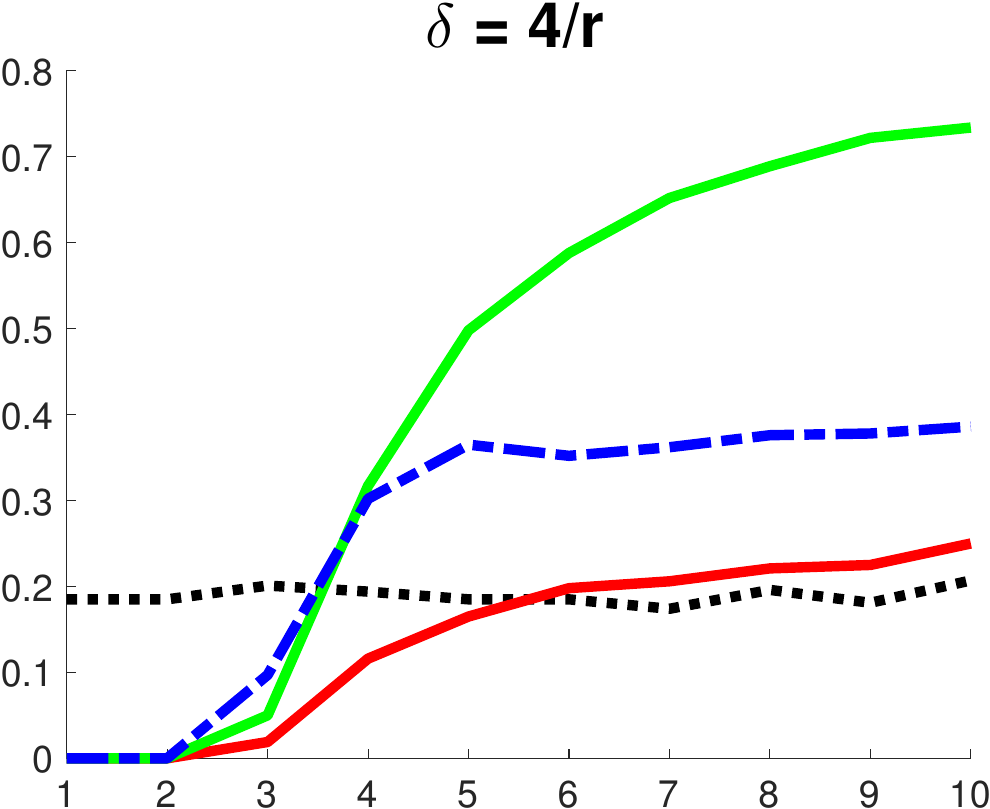}\\
\includegraphics[width=0.4\textwidth]{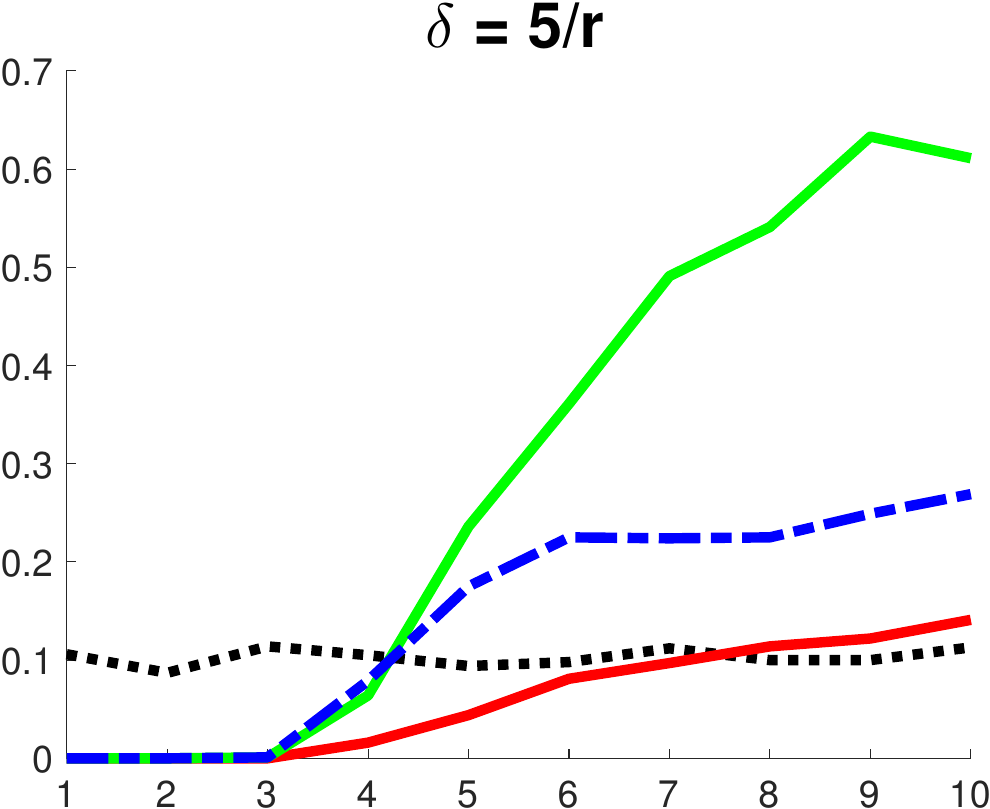}
\includegraphics[width=0.4 \textwidth]{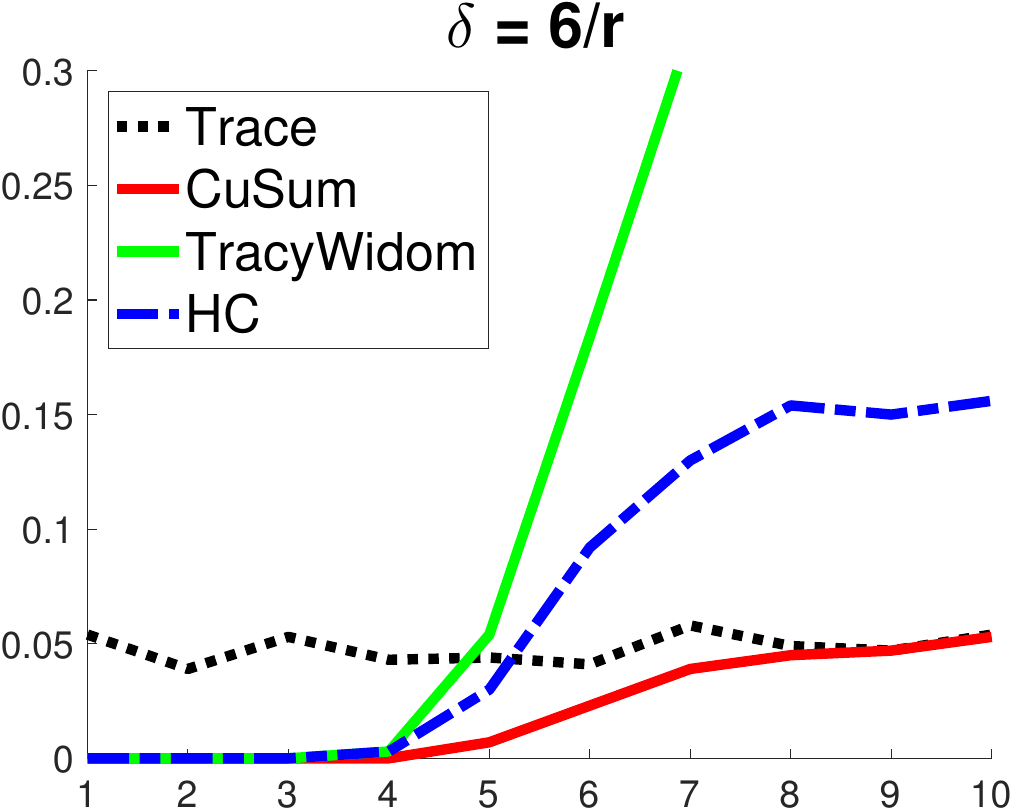}
\caption{Comparison of the ideal testing errors ($p=1200$, $n=1000$). The $x$-axis represents the number of spikes $r$, and for each $r=1,\cdots,10$, the strength of spikes is $\delta=c/r$, for a constant $c>0$.}   \label{fig:4tests}
\end{figure}

Third, we compare all four tests for a variety of $(r,\delta)$. Fixing a constant $c>0$, we consider $r\in\{1,\cdots,10\}$ and for each $r$ we let $\delta=c/r$. As $r$ changes, the performance of 
the trace test is roughly the same because 
\[
\frac{S_n - np}{\sqrt{2np}}   \mbox{$\approx N(0, 1)$  and  } \approx N(\frac{n(\delta r)}{\sqrt{2np}}, 1),\mbox{   under $H_0^{(n)}$ and  $H_1^{(n)}$, respectively}.  
\] 
So the trace test can be used as a benchmark. The results for $c\in \{3,4,5,6\}$ are displayed in Figure~\ref{fig:4tests}. When $r = 1$,  the TW test usually performs the best, and the trace test behaves 
unsatisfactorily.  
When $r$ is slightly larger, TW becomes less satisfactory and 
is inferior to all other three tests; in this setting, trace and CuSum
have the best performance.  

Recall that CuSum can be viewed as a hybrid of the trace and TW tests, we expect that  CuSum has both the advantage of TW and trace. This is confirmed by the 
numerical results: when $r$ is very small, CuSum behaves similarly to TW, and when $r$ is slightly larger, 
CuSum behaves similarly to trace.  Overall, it seems CuSum has the best performance, especially when the spikes are relatively strong (i.e., $c$ is large).

\begin{figure}[tb!] 
\centering
\includegraphics[width=0.42\textwidth]{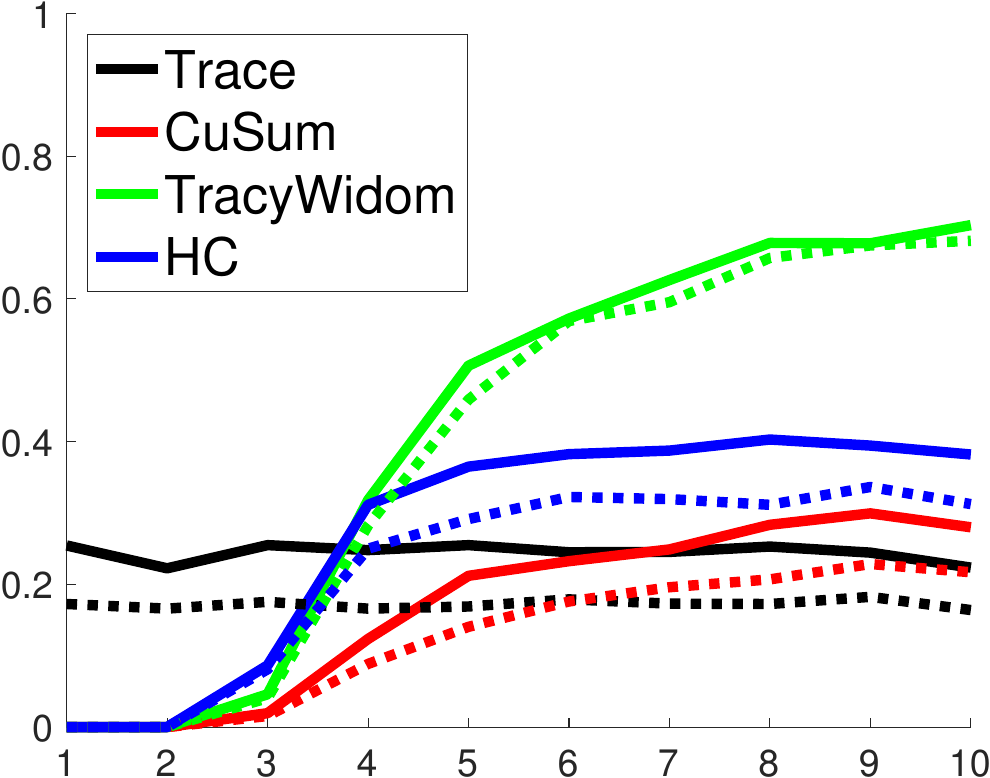}
\caption{The case where some non-spiked eigenvalues are smaller than $1$ under the alternative ($p=1500$, $n=800$). Dotted lines: all non-spiked eigenvalues are $1$. Solid lines: most non-spiked eigenvalues are $1$, except that $30$ of them are drawn uniformly from $[0.95,1]$.}  \label{fig:newSimu2}
\end{figure}

Fourth, we investigate a case where a small number of non-spiked eigenvalues are strictly smaller than $1$ under the alternative. Fix $(n,p)=(800,1500)$. For each of $r\in \{1,2,\ldots,10\}$, let $\delta = 5/r$. The null hypothesis is the same as before. For the alternative hypothesis, we generate $X_i\overset{iid}{\sim}N(0,\Sigma)$, where $\Sigma$ is a $p\times p$ diagonal matrix such that $\Sigma(j,j)=1+\delta$ for $1\leq j\leq r$, $\Sigma(j,j)$ is drawn uniformly from $[0.95,1]$ for $p-29\leq j\leq p$, and $\Sigma(j,j)=1$ for the remaining. The ideal testing errors of all four tests are displayed in Figure~\ref{fig:newSimu2} (solid lines). As a benchmark, we also consider a similar setting where the only difference is that $\Sigma(j,j)=1$ for $p-29\leq j\leq p$ under the alternative. The results for the bench mark setting are also in Figure~\ref{fig:newSimu2} (dotted lines). 

Compared with the benchmark setting, when some non-spiked eigenvalues are strictly smaller than $1$, the performance of all four tests become worse. The reason is that each empirical eigenvalue $\lambda_k$ is affected by all population eigenvalues (spiked and non-spiked ones). Decreasing of non-spiked eigenvalues will make $\lambda_k$'s smaller, so the powers of all four tests deteriorate. Such an effect is more mild for smaller $k$, which explains why the performance of TW test remains almost the same but the performance of trace test changes more significantly. By looking at the solid curves only, we can see that the CuSum test still has the best overall performance.

\begin{figure}[tb!] 
\centering
\includegraphics[width=0.328\textwidth]{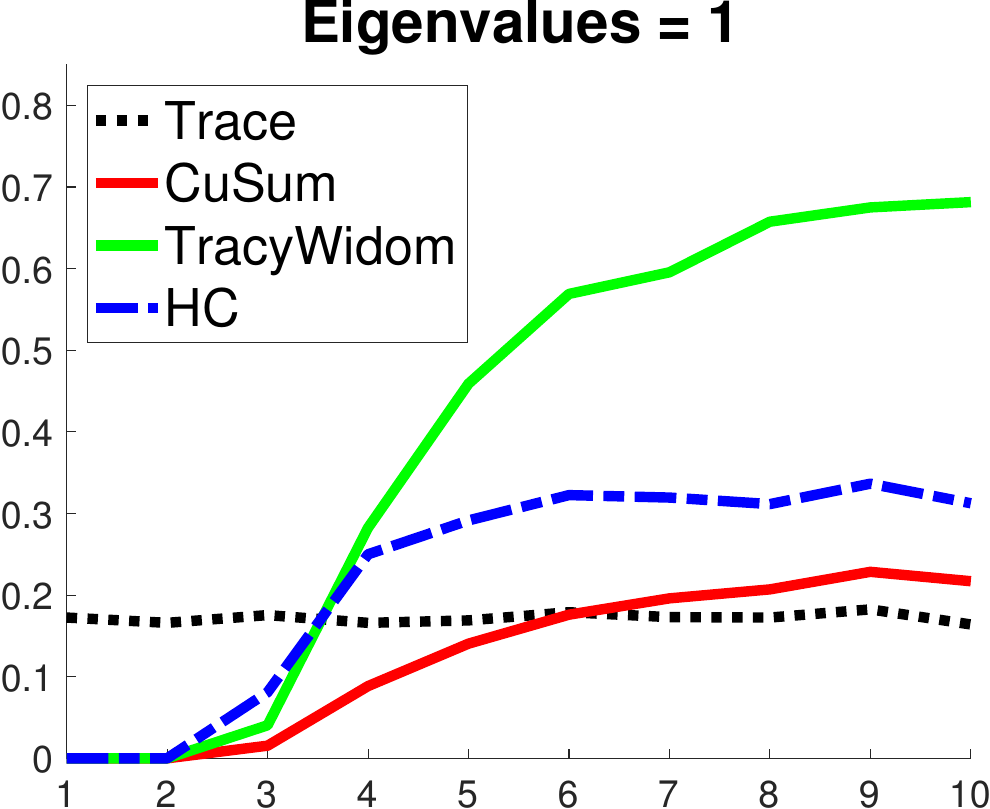}
\includegraphics[width=0.325\textwidth]{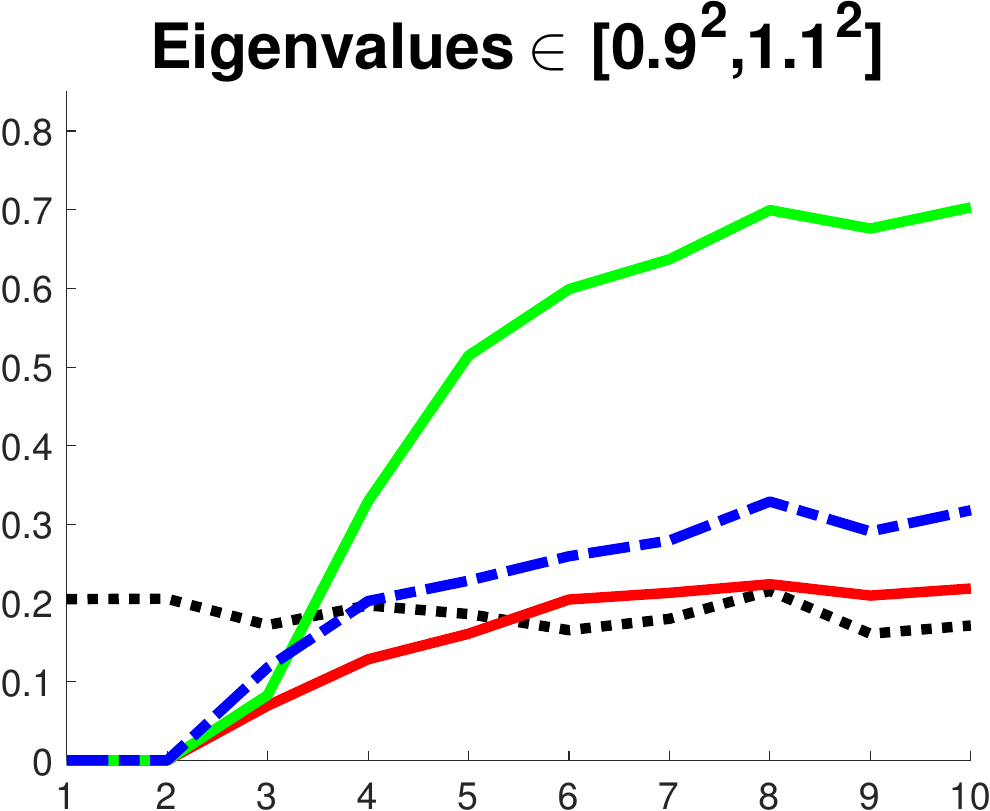}
\includegraphics[width=0.325\textwidth]{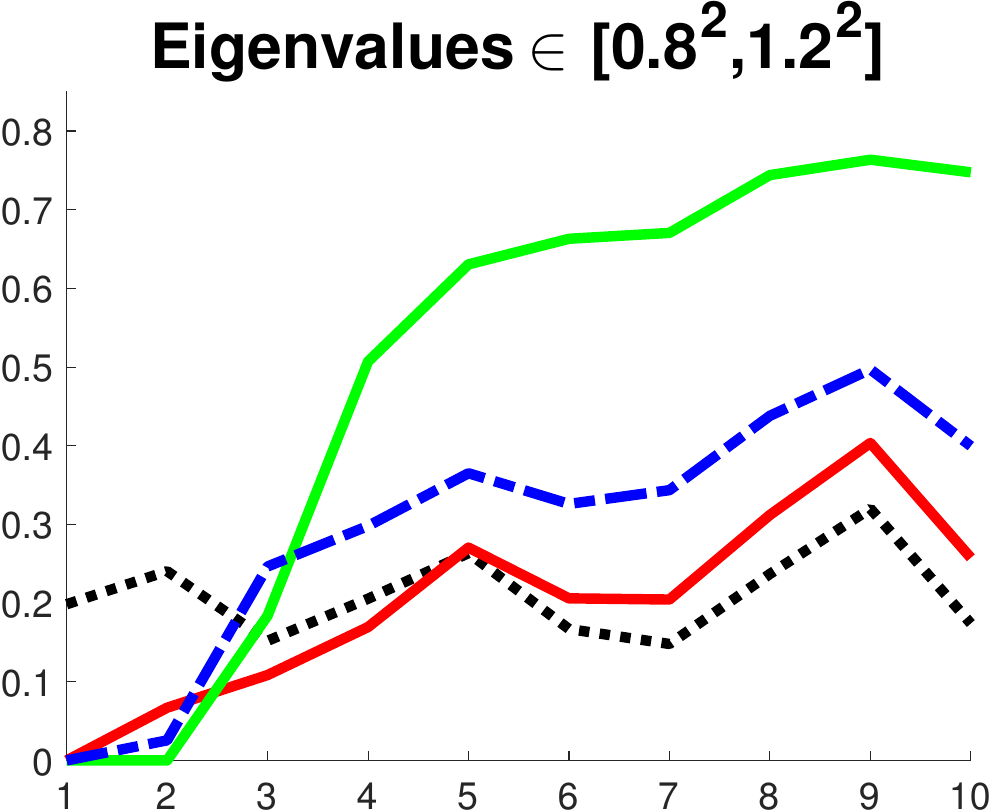}
\caption{The case where $\Sigma$ is not the identity matrix under the null ($p=1500$, $n=800$). The eigenvalues of $\Sigma$ are between $(1-b)^2$ and $(1+b)^2$, for $b=0$ (left), $0.1$ (middle) and $0.2$ (right).}  \label{fig:newSimu}
\end{figure}

Last, we consider the broader settings in Section~\ref{subsec:exten} where $\Sigma$ is not necessarily the identity matrix under the null. Fix $(n, p) = (800,1500)$. Given $b\in [0,1)$, we generate a $p\times p$ diagonal matrix $\Sigma=\diag(\sigma^2_1,\sigma_2^2,\ldots,\sigma^2_p)$ where $\sigma_j$'s are $iid$ drawn from the uniform distribution on $[1-b,1+b]$. Then, for each $(r,\delta)$, letting $\Sigma^*\in\mathbb{R}^{p,p}$ be the diagonal matrix with $\Sigma^*(j,j)=\Sigma(j,j)+\delta\cdot 1\{1\leq j\leq r\}$, we consider the testing problem where $X_i \stackrel{iid}{\sim} N(0, \Sigma)$ under the null and $X_i \stackrel{iid}{\sim} N(0, \Sigma^*)$ under the alternative. Although the null model is different from \eqref{model-null}, we still use the same implementation for all tests; for example, when computing the CuSum test statisitc, we still use $E_0(S_k)$ and $SD_0(S_k)$ simulated from model \eqref{model-null}. We consider $b\in \{0,0.1,0.2\}$, and for each choice of $b$, we run experiments for $r\in \{1,2,\ldots,10\}$ with $\delta=5/r$. The results are displayed in Figure~\ref{fig:newSimu}.

For $b=0.1$,  results are similar to those for $\Sigma=I_p$ (i.e., $b=0$). For $b=0.2$, when $r$ is small, results are similar to those for $\Sigma=I_p$; when $r$ gets larger, the performance of all tests becomes worse than that in the case of $\Sigma=I_p$ but is still reasonably good. Note that we haven't used any knowledge of $\Sigma$ in these tests. It is hopeful to further improve the performance by incorporating $\Sigma$, say, to compute the mean and standard deviation of $\lambda_k$ and $S_k$.

\section{Some results on Random Matrix Theory (RMT)} 
\label{sec:RMT} 
In this section, we present some Random Matrix Theory that are useful for our proofs. Some of them come from 
existing literature and some are newly developed. 

We write $Z\in \mathcal{G}(n,p)$ if $Z$ is an $n\times p$ random matrix with independent entries of $N(0,1)$. 
With probability $1$, the matrix $(1/n)Z'Z$ has $\min\{n,p\}$ distinct positive eigenvalues \cite{Uhlig94}, denoted as $\lambda_1> \lambda_2> \cdots>\lambda_{n\wedge p}>0$. 
We are interested in both cases of $p/n\goto 0$ and $p/n\goto \gamma >0$. 

Let $\mu^Z_n$ be the empirical spectral measure of $(1/n)Z'Z$: for any real value $x$, $\mu_n^Z((-\infty,x])=p^{-1}\sum_{j=1}^p1\{\lambda_j\leq x\}$, where $\lambda_1\geq \lambda_2\geq \cdots\geq\lambda_p$ are all the eigenvalues of $(1/n)Z'Z$. Let $\mu_{n,p}$ be the Marchenko-Pastur (MP) measure which has a point mass of $\max\{1-n/p, 0\}$ at $0$ and the probability density 
\beq  \label{MP}
\mu_{n,p}(x) = \frac{n}{2\pi x p} \sqrt{(x - a_{n,p}^{-})(a_{n,p}^+ - x)}, \qquad a_{n,p}^{-} < x < a_{n,p}^+, 
\eeq
where $a_{n,p}^{\mp}=(1\mp \sqrt{p/n})^2$. 
For any $1\leq k\leq (n\wedge p)$, let $q_k=q_k^{(n,p)}$ be such that 
\beq  \label{Definequantile}
\int_{q_k}^{a_{n,p}^+} \mu_{n,p}(x)dx = k/n. 
\eeq
We also define $\mu_{\gamma}$, the Marchenko-Pastur (MP) measure associated with $\gamma$, which has a point mass of $\max\{1-1/\gamma, 0\}$ at $0$ and the probability density 
\beq\label{MP2}
\mu_\gamma(x) = \frac{1}{2\pi x\gamma} \sqrt{(x - a_{\gamma}^{-})(a_{\gamma}^+ - x)}, \qquad a_{\gamma}^{-} < x < a_{\gamma}^+,
\eeq
where $a_\gamma^{\mp}=(1\mp\sqrt{\gamma})^2$. Note that $\mu_{n,p}$ is a special case of $\mu_{\gamma}$ with $\gamma=p/n$.

\subsection{Deviation of linear eigenvalue statistics}
In this section, we assume $p\leq n$ without loss of generality. 
For any function $f:\mathbb{R}\to \mathbb{R}$, we call
\[
\mu^Z_n(f)= \frac{1}{p}\sum_{i=1}^{p} f(\lambda_i)
\] 
the {\it linear eigenvalue statistic} associated with $f$. The following theorem establishes a deviation inequality for such linear statistics when $f'$ is uniformly bounded. It is proved in the appendix. 
\begin{thm} \label{thm:deviation} 
Suppose $Z\in\mathcal{G}(n,p)$, $p\leq n$ and that $f$ is continuously differentiable satisfying $\|f'\|_\infty\equiv \sup_{x}|f'(x)|<\infty$. For any $t>0$, 
\[
P\left( |\mu_n^Z(f) - E[\mu^Z_n(f)]| > t \|f'\|_\infty/\sqrt{np} \right) \leq 2\inf_{y\geq\sqrt{2\log(2)}}\bigl\{e^{-\frac{t^2}{c_0y^2}} + 2e^{-y^2}\bigr\},
\]
where $c_0>0$ is an absolute constant. 
\end{thm} 

A similar large-deviation inequality is given in \cite[Corollary 1.8(b)]{GZ00}, but Theorem~\ref{thm:deviation} improves it in two folds. First, Theorem~\ref{thm:deviation} gives a sharper rate in the case $p/n\goto 0$. Second, Theorem~\ref{thm:deviation} only requires that $f$ is Lipschitz, while \cite[Corollary 1.8(b)]{GZ00} requires that $f(x^2)$ is Lipschitz. However, we pay a price of the extra term $2e^{-y^2}$ on the right hand side. 

A direct application of Theorem~\ref{thm:deviation} is to derive a large deviation result for $\sum_{i=1}^p(\lambda_i-1)^m$ in the case $p/n\goto 0$, where $m$ is a positive integer. The following lemma is proved in the appendix. 

\begin{lemma}  \label{lem:eigen-moment}
Suppose $Z\in\mathcal{G}(n,p)$ and $p/n\goto 0$ as $n\goto\infty$. For any fixed integer $m\geq 1$, as $p\goto\infty$, with probability at least $1-o(n^{-1})$, 
\[
\Bigl| \sum_{i=1}^p \{(\lambda_i-1)^m- E[(\lambda_i-1)^m]\}\Bigr|  \leq C\log(n) \left[\frac{p\vee \log(n)}{n}\right]^{m/2}. 
\]
\end{lemma}

Although large-deviation inequalities for $\sum_{i=1}^p (\lambda_i-1)^m$ could be obtained from existing results on polynomials of independent normals (e.g., \cite{Wolff, Vu}), those inequalities are expressed in terms of tensor norms or subgraph counts, from which it is hard to get the explicit rate. We obtain Lemma~\ref{lem:eigen-moment} almost for free by applying Theorem~\ref{thm:deviation} to  $f(x)=(x-1)^m$ (with a proper truncation).  

\subsection{Concentration of empirical spectral measure}
In this section, we focus on the case that $p/n\goto\gamma$, where $\gamma>0$. 
It is well known that $\mu^Z_{n}\goto\mu_{\gamma}$ weakly. 
The following theorem, which is a special case of \cite[Theorem 1.1]{KSdist}, gives its convergence rate in terms of the Kolmogorov-Smirnov distance between distribution functions of $\mu^Z_{n}$ and $\mu_{n,p}$. 
\begin{thm}\label{thm:KS}
Suppose $Z\in\mathcal{G}(n,p)$ and $p/n\goto \gamma$ for a constant $\gamma>0$. Let $F_n^Z$ and $F_{n,p}$ be the CDF of $\mu_n^Z$ and $\mu_{n,p}$, respectively. For any $\zeta>0$, there exist positive constants $C$ and $c$ which depend on $(\zeta,\gamma)$, such that   
\[
P\left( \sup_{x\in \mathbb{R}}|F^Z_n(x) - F_{n,p}(x)| > n^{-1}\omega_{n,\zeta}^6 \right) \leq C_\zeta \exp( -c_\zeta  \omega_{n,\zeta}), 
\]
where $\omega_{n,\zeta}=\log(n)[\log(\log(n))]^\zeta$.
\end{thm} 

The next lemma gives an explicit expression of $\int \log(t-x)\mu_\gamma(x)dx$, for $t>a_\gamma^+$, and is proved in the appendix. We have not found such a result in the literature (although a formula for $\int \log(t+x)\mu_\gamma(x)dx$ is given in \cite{Verdu}). 
\begin{lemma} \label{lem:MPintegral}
For any $\gamma>0$ and $t>(1+\sqrt{\gamma})^2$,
\begin{align}
\int &\log(t-x)\mu_\gamma(x) dx = \frac{\gamma-1}{2\gamma} \log\big((\gamma+1)A_t - (\gamma-1)B_t+4\gamma\big)\cr
& + \frac{A_t- B_t}{2\gamma} -\frac{\gamma +1}{2\gamma}\log(A_t- B_t) +\frac{1}{\gamma}\log(2) + \frac{\gamma+1}{2\gamma}\log(\gamma), 
\end{align}
where $A_t = t-\gamma -1$ and  $B_t = \sqrt{(t-\gamma-1)^2-4\gamma}$.
\end{lemma}

\subsection{Behavior of eigenvalues and their cumulative sums}
In this section, we study $\lambda_k$ and $S_k=\sum_{j=1}^k\lambda_j$. Again,  we focus on the case that $p/n\goto \gamma >0$. Consider $\lambda_k$ first. The following theorem is a special case of \cite[Theorem 3.3]{PillaiYin}.

\begin{thm}  \label{thm:eig1}
Suppose $Z\in\mathcal{G}(n,p)$ and $p/n\goto \gamma$ for a constant $\gamma>0$ and $\gamma\neq 1$.
Write $L_n=[\log(n)]^{\log(\log(n))}$. For any $\zeta>0$, there exists a constant $C_\zeta>0$ such that 
\[
P\left( \cup_{1\leq k\leq n}\left\{  |\lambda_k -q_k|> L^{C_\zeta}_n n^{-2/3}[k\wedge (n+1-k)]^{-1/3}\right\}  \right)\leq n^{C_\zeta} \exp( -L_n^{\zeta}). 
\]
\end{thm}

In the case $\gamma=1$, Theorem~\ref{thm:eig1} continues to hold, except that $1\leq k\leq n$ is now replaced with $1\leq k\leq (1-\eta)n$, for a positive constant $\eta\in (0,1)$; see the remarks in the end of \cite[Section 4]{PillaiYin}. We note that a direct corollary of Theorem~\ref{thm:KS} gives similar results but only for $\omega_{n,\zeta}^6\leq k\leq n-\omega_{n,\zeta}^6$, however, Theorem~\ref{thm:eig1} holds for the whole range $1\leq k\leq n$.   

As a corollary of Theorem~\ref{thm:eig1}, 
\beq
|E(\lambda_k)-q_k|\leq L^{c}_n n^{-2/3}\tilde{k}^{-1/3}, \quad SD(\lambda_k)\leq L^{c}_n n^{-2/3}\tilde{k}^{-1/3}, 
\eeq
where $\tilde{k}=k\wedge(n+1-k)$ and $c>0$ is a properly large constant. The null behavior of $HC_n^+$ follows immediately. Now, we consider $HC_n^*$. To access its null behavior,  we additionally need a lower bound for $SD(\lambda_k)$. 
Unfortunately, we have not found a desirable lower bound of  $SD(\lambda_k)$ in the existing literature. According to \cite{Sandrine2013, Su}, we conjecture that
\[
|E(\lambda_k)-q_k|\asymp n^{-2/3}\tilde{k}^{-1/3}, \quad SD(\lambda_k)\asymp \sqrt{\log(n)} n^{-2/3}\tilde{k}^{-1/3}. 
\]

The following theorem gives a similar result for the cumulative sum $S_k=\sum_{j=1}^k\lambda_j$. It is a generalization of the result in \cite{Bao-PS} and proved in \cite{Bao}. 
\begin{thm}  \label{thm:eig2}
Suppose $Z\in\mathcal{G}(n,p)$ and $p/n\goto \gamma$ for a constant $\gamma>0$.
Write $L_n=[\log(n)]^{\log(\log(n))}$. For any $\zeta>0$, there exists a constant $C_\zeta>0$ such that 
\[
P\Bigl( \cup_{1\leq k\leq n}\Bigl\{  |S_k -\sum_{i=1}^k q_i|> L^{C_\zeta}_n (k/n)^{2/3} \Bigr\}  \Bigr)\leq n^{C_\zeta} \exp( -L_n^{\zeta}). 
\]
\end{thm}
The discussions of the null behavior of $CS_n^+$ and $CS_n^*$ are similar to those of $HC_n^+$ and $HC_n^*$, so we omit them. 
\section{Proof of main theorems} \label{sec:main} 
We prove the main results in Section~\ref{sec:intro}. 
For Theorem~\ref{thm:CS}, 
the first claim is a direct result of Theorem~\ref{thm:eig2} (so the level of the test tends to $0$) and that $CS_n^+\geq S_n-E_0[S_n]$ (so the power of the test is no smaller than the power of the trace test), and the second claim follows from the lower bound argument and the null behavior of eigenvalues (Theorem~\ref{thm:eig2}). Similarly, Theorem~\ref{thm:HC} follows from Theorem~\ref{thm:eig1} and the lower bound argument. Theorems~\ref{thm:LeCam2}-\ref{thm:LeCam3} are direct results of Lemmas~\ref{lemma:LeCam2}-\ref{lemma:LeCam3}, respectively, and Theorem~\ref{thm:covert} follows directly from Theorem~\ref{thm:phase}. 

What remains is to prove Theorem~\ref{thm:phase}, Corollary~\ref{cor:TW}, Lemmas~\ref{lemma:LeCam2}-\ref{lemma:LeCam3}, and Theorem~\ref{thm:critical1}.

\subsection{Preliminary I: Gaussian proxy model and its intimacy to RWS}
At the heart of our analysis is to consider a proxy testing problem where we replace 
the alternative hypothesis $H_1^{(n)}$ by   
\begin{equation} \label{model-tildealt}  
\widetilde{H}_1^{(n)}:   \qquad  \widetilde{\Sigma} =  [I_p - (\delta/p)  Y Y']^{-1},  \;\;\;\;  Y = Z \cdot 1\{\|Z\| \leq \frac{1}{2}\sqrt{p/\delta}\},  \footnote{Requiring that $\|Z\| \leq  \frac{1}{2}  \sqrt{p/\delta}$ seems to be a nuisance but is necessary: otherwise the matrix $[I_p - (\delta/p) YY']$ may not be invertible; note that $P\big(\|Z\| > (1/2) \sqrt{p/\delta}\big)$ is negligible.}
\end{equation} 
where $Z \in \mathbb{R}^{p,r}$ is the matrix that has $iid$ $N(0,1)$ entries.    
\bed  
We call (\ref{model-cal2}) and (\ref{model-tildealt})  the Gaussian proxy model of RWS. 
\eed 

The next lemma says that two models are close to each other, for a wide region in the parameter space. Let $LR_n$ and $\widetilde{LR}_n$ be the likelihood ratios associated with the RWS and the proxy model, respectively. Let $E_0$ stand for the expectation under $H_0^{(n)}$.  
\begin{lemma} \label{lem:differ} 
({\it Intimacy of two models}).  Consider the testing problem \eqref{model-null}-\eqref{model-alt} where \eqref{model-cal2} holds.   If $\alpha>\max\{1-5\beta/4, (1-\beta)/2\}$, then as $n \goto \infty$, 
$E_0|LR_n-1| - E_0|\widetilde{LR}_n-1| \goto 0$.
\end{lemma} 
See Figure \ref{fig:phase} for an illustration.

This lemma is proved in Section~\ref{subsec:proof-differ}. The proof uses Le Cam's ``comparison of experiments", an idea explained in Section~\ref{subsec:exten}. We conjecture a stronger result of $E_0|LR_n - \widetilde{LR}_n| \goto 0$ holds; see Section~\ref{subsec:proof-differ} for more discussions.  

Our motivation of introducing the proxy model is that it is much easier to analyze than the RWS. 
Let $F(Q)$ and $F(Y)$ be the CDF of $Q$ in RWS and $Y$ in proxy model.  
It will be shown in Section~\ref{subsec:proxy} that  
\begin{align} 
&LR_n = (1 - \delta)^{\frac{nr}{2}}  \int   \mexp(\frac{\delta}{2} \tr(\hat{\Sigma} QQ')) d F(Q),  \label{DefineLR} \\ 
&\widetilde{LR}_n \approx (1 - \delta)^{\frac{nr}{2}}  e^{\frac{n r \delta}{2(1 - \delta)} (1 - \frac{r \delta}{p(1 - \delta)})}  \int e^{\frac{\delta n}{2p} \tr((\hat{\Sigma} - \frac{1}{1 - \delta} I_p)  YY')} d F(Y).
\label{tildeRadd} 
\end{align} 
The main technical hurdle we face in RWS is that, for the integral in (\ref{DefineLR}), it is hard to integrate $Q$ out. 
However, such a hurdle is removed by using a proxy model:  for the integral in (\ref{tildeRadd}), 
we can 
integrate $Y$ out and derive a (relatively) simple formula, because the entries of $Y$ are (approximately) $iid$ $N(0,1)$ and the integral can be expressed using linear spectral statistics of a standard Wishart matrix.


\subsection{Preliminary II: Likelihood ratio of the proxy model} \label{subsec:proxy}
We justify the approximation in \eqref{tildeRadd} and calculate the integral over $Y$. The justification of \eqref{tildeRadd} is relatively short and conventional, however, for $LR_n$, it is unclear how to derive such an approximation.

Let $Y \in \mathbb{R}^{p,r}$ be a random matrix of {\it iid} $N(0,1)$ entries. Since $r\ll p$, with probability $1$, $p^{-1}YY'$ has $r$ distinct positive eigenvalues \cite{Uhlig94}, denoted as 
$\eta_1 > \eta_2 \ldots  >\eta_r>0$. Then, the proxy model \eqref{model-tildealt} is equivalent to  
\begin{equation} \label{md-tildealt}  
\widetilde{H}_1^{(n)}: \quad \widetilde{\Sigma} =  \begin{cases}
[I_p - \delta (p^{-1} Y Y')]^{-1}, & \mbox{if } \delta\eta_1\leq 1/2,\\
I_p, & \mbox{otherwise}. 
\end{cases}
\end{equation} 
Let  $F(Y)$ be the joint distribution function of $Y$. 
The likelihood ratio $\widetilde{LR}_n$ associated with \eqref{md-tildealt} is 
\begin{align} \label{DefinetildeLR} 
\widetilde{LR}_n   
&= (1-\delta)^{\frac{nr}{2}} \int_{\delta\eta_1\leq 1/2} \mexp\left\{ \frac{\delta n}{2p}\tr(\hat{\Sigma}YY')  
+ \frac{n}{2}g_n(Y,\delta,r) \right\} d F(Y) \cr
&+P(\delta\eta_1> 1/2), 
\end{align}
where $g_n(Y, \delta, r) = \sum_{k = 1}^r  \log(1 - \frac{\delta}{1 - \delta}(\eta_k -1))$.  

We then show that $\widetilde{LR}_n$ can be approximated by $\widetilde{R}_n$, to be introduced below, which has a relatively simple form. Introduce $g_n^*(Y, \delta, r)  = - \frac{\delta}{1 - \delta} \sum_{k=1}^r (\eta_k -1) -  \frac{\delta^2}{2 (1 - \delta)^2} \sum_{k = 1}^r E[(\eta_k -1)^2]$; it is the Taylor expansion of $g_n(Y, \delta, r)$ up to the second term, except for $(\eta_k-1)^2$ is replaced by $E[(\eta_k - 1)^2]$. 
Substituting $g_n$ with $g_n^*$ gives a proxy of $\widetilde{LR}_n$:\footnote{We also remove $\delta\eta_1<1$ in the integral of \eqref{DefinetildeR} so that it has a closed form.}
\begin{align}  \label{DefinetildeR}
\widetilde{R}_n &= 
(1 - \delta)^{\frac{nr}{2}} 
\int \mexp\left\{\frac{\delta n}{2p}\tr(\hat{\Sigma}YY')  
+ \frac{n}{2} g_n^*(Y, \delta, r)  \right\} d F(Y)\cr
&= (1 - \delta)^{\frac{nr}{2}} b_n(\delta, r)
\int  \mexp\left\{ \frac{\delta n}{2p}\tr\bigl((\hat{\Sigma}-\frac{1}{1-\delta}I_p)YY'\bigr)\right\} dF(Y),
\end{align}
where $b_n(\delta,r)=e^{ \frac{nr\delta}{2(1-\delta)} - \frac{nr^2\delta^2}{4p(1-\delta)^2} }$. Here the second equality follows from that $\sum_{k=1}^rE[(\mu_k-1)^2]=r^2/p$ \cite[Lemma 2.9]{Verdu}.
Noting that the above $\widetilde{R}_n$ is finite only when $\frac{\delta n}{p}(\lambda_1-\frac{1}{1-\delta})<1$, 
we simply set $\widetilde{R}_n=1$ when $\frac{\delta 
n}{p}(\lambda_1-\frac{1}{1-\delta})\geq 1$.

The following lemma is proved in the appendix. 
\begin{lemma} \label{lem:approx}  
Consider the testing problem \eqref{model-null} and \eqref{model-tildealt} with parameters as in \eqref{model-cal2}.  
If $\alpha>\max\{(4-5\beta)/6, (1-\beta)/2\}$, then 
$E_0|\widetilde{LR}_n - \widetilde{R}_n|\leq Cn^{1-\min\{2\alpha+\beta,\; 3\alpha +5\beta/2-1\}}$. 
\end{lemma}

A nice property of $\widetilde{R}_n$ is that we can integrate out $Y$ explicitly in \eqref{DefinetildeR}. Denote by $\mu^X_n$ the empirical spectral measure of $\hat{\Sigma}$ under the null, i.e., $\mu^X_n=p^{-1}\sum_{j=1}^p \delta_{\lambda_j}$, where $\delta_a$ represents the point mass at $a$ and $\lambda_1\geq \lambda_2\cdots\geq \lambda_p$ are all eigenvalues of $\hat{\Sigma}$ (including zero ones). 
Let $\mu_{n,p}$ be the Marchenko-Pastur (MP) measure with parameter $\gamma=p/n$ as in \eqref{MP}. 
The next lemma gives two equivalent expressions of $\widetilde{R}_n$. 

\begin{lemma} \label{lem:tildeR}
Consider the testing problem \eqref{model-null} and \eqref{model-tildealt} with parameters as in \eqref{model-cal2}. Let $\psi_n(\lambda)=\frac{p}{n\delta}\log(1- \frac{\delta n}{p}(\lambda-\frac{1}{1-\delta}))$. When $\frac{\delta n}{p}(\lambda_1-\frac{1}{1-\delta})<1$, 
\begin{eqnarray} 
&& \widetilde{R}_n = (1-\delta)^{\frac{nr}{2}} b_n(\delta, r) 
\mexp\Big\{-\frac{r}{2}\sum_{j=1}^p \log\bigl( 1-\frac{\delta n}{p}(\lambda_j-\frac{1}{1-\delta}) \bigr)\Big\} \label{proxy1a}\\
&&= e^{-\frac{nr^2\delta^2}{4p(1-\delta)^2}} \mexp\left\{-  \frac{nr\delta}{2} \left[\int \psi_n(\lambda) \mu_n^X(d\lambda) - \int \psi_n(\lambda) \mu_{n,p}(d\lambda)\right]\right\}, \label{proxy1b}
\end{eqnarray} 
where $b_n(\delta,r)$ is the same as in \eqref{DefinetildeR}. 
\end{lemma}

We have seen that
\[
LR_n\approx \widetilde{LR}_n\approx \widetilde{R}_n,   
\]
and according to \eqref{proxy1b}, the key of analyzing $\widetilde{R}_n$ is to characterize the convergence of empirical spectral measure $\mu_n^X$ to the Marchenko-Pastur measure $\mu_{n,p}$. This is well-studied in Random Matrix Theory.

\subsection{Proof of Theorem~\ref{thm:phase}}   \label{subsec:proof-LB}
We only prove the first claim. By Neyman-Pearson lemma, it suffices to show that when $\alpha+\beta> 1$, $E_0|LR_n -1|\goto 0$. Recall that $\widetilde{LR}_n$ is the likelihood ratio associated with the Gaussian proxy model, and $\widetilde{R}_n$ is a proxy of $\widetilde{LR}_n$ introduced in Section~\ref{subsec:proxy}. By Lemma~\ref{lem:differ} and Lemma~\ref{lem:approx}, 
$E_0|LR_n-1|= E_0|\widetilde{R}_n-1|+o(1)$. So it suffices to show that 
\beq \label{thm-phase-0}
E_0|\widetilde{R}_n-1|\goto 0. 
\eeq

We now show \eqref{thm-phase-0}. Write $\psi_n(\lambda)=\frac{p}{n\delta}\log(1- \frac{\delta n}{p}(\lambda-\frac{1}{1-\delta}))$. Let $\mu_n^X$ be the empirical spectral measure associated with $\hat{\Sigma}$, and let $\mu_{n,p}$ be the Marchenko-Pastur (MP) measure as in Section~\ref{sec:RMT}.
Using \eqref{proxy1b} of Lemma~\ref{lem:tildeR},  we have
\beq \label{thm-phase-3}
\log(\widetilde{R}_n) = -  \frac{nr\delta}{2} \left[\int \psi_n(\lambda) \mu_n^X(d\lambda) - \int \psi_n(\lambda) \mu_{n,p}(d\lambda)\right] -\frac{nr^2\delta^2}{4p(1-\delta)^2}, 
\eeq
Let $F_n^X$ and $F_{n,p}$ be the CDF associated with $\mu_n^X$ and $\mu_{n,p}$, respectively. Introduce $E_n=\{X\in\mathbb{R}^{n,p}: \sup_{x\in\mathbb{R}}|F_n^X(x)-F_{n,p}(x)|\leq n^{-1}\log^6(n)\log^6(\log(n))\}$, and let $I_{E_n}$ be the indicator that $X\in E_n$, similarly for $I_{E_n^c}$. By Theorem~\ref{thm:KS}, as $n\goto\infty$, we have $P_0(X\in E_n^c)=O(e^{-C\log(n)\log(\log(n))})=o(n^{1-\alpha-\beta})$. Since $E_0(\widetilde{R}_n)\leq E_0(\widetilde{LR}_n)+o(1)$, which is bounded, we have $E_0(\widetilde{R}_n\cdot I_{E_n^c})=o(1)$. As a result, 
\beq \label{thm-phase-1}
E_0(| \widetilde{R}_n - 1|\cdot I_{E_n^c}) \leq E_0 (\widetilde{R}_n\cdot I_{E_n^c})
+P_0(X\in E_n^c)\goto 0. 
\eeq
Below, we consider $E_0(|\widetilde{R}_n-1|\cdot I_{E_n})$. 
Over the event $E_n$, both the supports of $\mu_n^X$ and $\mu_{n,p}$ are strictly contained in $A_\eps=[(1-\sqrt{\gamma})^2-\eps, (1+\sqrt{\gamma})^2+\eps]$, for any $\eps\in (0,1)$. Furthermore, 
$\psi_n'(\lambda)=-[1-\frac{\delta n}{p}(\lambda- \frac{1}{1-\delta})]^{-1}$; it follows that $\sup_{\lambda\in A_\eps}|\psi'_n(\lambda)|\leq 2$. Combining the above results, for $X\in E_n$, 
\begin{align*}
\left|\int \psi_n(\lambda) \mu_n^X(d\lambda) - \right.&\left. \int \psi_n(\lambda) \mu_{n,p}(d\lambda)\right|\leq \sup_{\lambda\in A_\eps}|\psi'_n(\lambda)|\cdot \sup_{\lambda} |F_n^X(\lambda)-F_{n,p}(\lambda)|\cr
&\leq 2n^{-1}\log^6(n)\log^6(\log(n)). 
\end{align*}
Plugging it into \eqref{thm-phase-3} gives $|\log(\widetilde{R}_n)|\leq L_n\delta r+C\delta^2r^2\leq L_n\delta r$, where $L_n$ is a multi-$\log(n)$ term. This implies that for $X\in E_n$, 
$|\widetilde{R}_n-1|\leq \max\{|e^{L_n\delta r}-1|, |e^{-L_n\delta r}-1|\}\leq L_n\delta r$. 
As a result, 
\beq \label{thm-phase-2}
E_0(|\widetilde{R}_n-1|\cdot I_{E_n})\leq L_n\delta r =  L_n n^{1-\alpha-\beta}. 
\eeq
Combining \eqref{thm-phase-1}-\eqref{thm-phase-2} gives \eqref{thm-phase-0}. 
\qed

\subsection{Proof of Corollary~\ref{cor:TW}} 
The first claim follows immediately from the null behavior of $\lambda_1$. We only prove the second claim, that is, if $\alpha>2/3$, the power of the TW test tends to $0$. 

We first show that, if we fix all other parameters in the RWS \eqref{model-alt} but increase $r$, the power of the TW test always increases. For $r_1, r_2$ such that $r_2>r_1$, let $Q_1$ and $
Q_2$ be uniformly distributed over $\mathbb{S}(p,r_1)$ and $\mathbb{S}(p,r_2)$, respectively, and write $\Sigma_1=I_p +\frac{\delta}{1-\delta} Q_1Q_1'$ and $\Sigma_2=I_p +\frac{\delta}{1-\delta} Q_2Q_2'$. Consider the two RWS models $X_i|Q_1\overset{iid}{\sim}N(0, \Sigma_1)$ and $X_i|Q_2\overset{iid}{\sim}N(0, \Sigma_2)$. 
Write $Q_2=[\tilde{Q}_2, R]$, where $\tilde{Q}_2$ is the submatrix of $Q_2$ formed by its first $r_1$ columns. It is known that $\tilde{Q}_2$ is uniformly distributed over $\mathbb{S}(p,r_1)$. As a result, 
\[
\Sigma_1 \overset{(d)}{=} I_p + \frac{\delta}{1-\delta} \tilde{Q}_2\tilde{Q}_2', \qquad \Sigma_2 = I_p + \frac{\delta}{1-\delta} (\tilde{Q}_2\tilde{Q}_2' +RR'). 
\]
For each realization of $Q_2=[\tilde{Q}_2, R]$, $\Sigma_2-\Sigma_1$ is positive semi-definite. So, for any value $a>0$, the probability that $\lambda_1>a$ under the first RWS model is no larger than the probability that $\lambda_1>a$ under the second RWS model. This implies that the power of the TW test is larger in the second model.

Then, it suffices to consider the case that $r=n$. In this case, $\Sigma=\frac{1}{1-\delta} I_p$ under the alternative. By \cite{Johnstone2001}, if $X_i\overset{iid}{\sim}N(0, I_p)$, then $\frac{n^{2/3}(\lambda_1 - a_{n,p}^+)}{\sigma_{n,p}}$ converges weakly to the Tracy-Widom distribution, where $a_{n,p}^+=(1+\sqrt{p/n})^2$ and $\sigma_{n,p} = (1+\sqrt{p/n})(1+\sqrt{n/p})^{1/3}$, both converging to a constant as $n\goto\infty$. Therefore, when $\Sigma=\frac{1}{1-\delta} I_p$,  
\[
W\equiv \frac{n^{2/3}[(1-\delta) \lambda_1 - a_{n,p}^+]}{\sigma_{n,p}} \mbox{ converges weakly to the TW law}.  
\]
The rejection region is that $\lambda_1> a_{n,p}^+ + n^{-2/3}\sigma_{n,p}[3\log(n)]^{2/3}$. Hence, under the alternative, the rejection probability is equal to 
\[
P\left( W> - (n^{2/3}\delta) \cdot \sigma^{-1}_{n,p} a_{n,p}^+ + (1-\delta)[3\log(n)]^{2/3} \right). 
\]
Since $\alpha>2/3$, the term $n^{2/3}\delta=n^{2/3-\alpha}$ is dominated by the term $[\log(n)]^{2/3}$, so the above probability tends to $0$. \qed

\subsection{Proof of Lemma~\ref{lemma:LeCam2}}
Let $A=\widetilde{\Sigma}-\Sigma$ and $Z_i\overset{iid}{\sim}N(0, A)$, $1\leq i\leq n$, independent of $X_i$'s. Write $X=[X_1,\cdots, X_n]'$ and $Z=[Z_1,\cdots, Z_n]'$. It is seen that for $s\in\{0,1\}$, if $X$ follows the model in $H^{(n)}_s$, then $\widetilde{X}=X+Z$ follows the model in $\widetilde{H}_s^{(n)}$. Therefore, letting $h(Z)$ be the joint density of $Z$, we have 
\[
\tilde{f}_s(\widetilde{X}) = \int f_s(\widetilde{X}-Z)h(Z)dZ, \qquad s=0, 1. 
\]
It follows that
\begin{align*}
\|\tilde{f}_0 - \tilde{f}_1 \| &= \int \Big|\int f_0(\widetilde{X}-Z)h(Z)dZ - \int f_0(\widetilde{X}-Z)h(Z)dZ \Big|d\widetilde{X}
\cr 
&\leq \int\int |f_0(\widetilde{X}-Z)-f_1(\widetilde{X}-Z)|d\widetilde{X} h(Z)dZ\cr
&= \int \|f_0-f_1\| h(Z)dZ = \|f_0 - f_1\|, 
\end{align*}
where the third line is due to a variable change $X=\widetilde{X}-Z$. 
\qed

\subsection{Proof of Lemma~\ref{lemma:LeCam3}}

We prove the following lemma, and Lemma~\ref{lemma:LeCam3} is a special case of it. 
\begin{lemma} \label{lem:monotone2} 
({\it Monotonicity of $L^1$-distance}).   
Let $\Sigma_0(Q)$, $\Sigma_1(Q)$, and $\Sigma_2(Q)$ be three covariance matrix models   
indexed by a random quantity $Q$ such that 
\begin{itemize}
\item For any non-random orthogonal matrix $U$, $U'\Sigma_\ell(Q)U$ follows the same distribution as $\Sigma_{\ell}(Q)$, $\ell =0,1,2$.  
\item Given any realization $Q$, $\Sigma_1(Q)\preceq  \Sigma_2(Q)$. 
\item Given any two realizations $Q$ and $\widetilde{Q}$, $\Sigma_0(Q)\preceq  \Sigma_1(\widetilde{Q})$. 
\end{itemize}    
For each $\ell \in \{0, 1, 2\}$, consider a model where $(X_i | Q)$ are iid samples from $N(0, \Sigma_{\ell}(Q))$, $1 \leq i \leq n$, and let $f_{\ell}$ be the joint density of $\{X_i\}_{i = 1}^n$. We have $\|f_1 - f_0\|_1 \leq \|f_2 - f_0\|_1$. 
\end{lemma}
{\it Proof.} 
Consider the problem of testing $H_1: X_i|Q\overset{iid}{\sim}N(0, \Sigma_1(Q))$ against $H_0: X_i|Q\overset{iid}{\sim}N(0, \Sigma_0(Q))$. 
By Neyman-Pearson lemma, $1-(1/2)\|f_1 - f_0\|_1$ is a lower bound for the sum of type I and type II errors of any test, and this lower bound is achieved by the likelihood ratio test which rejects $H_0$ if and only if
\beq \label{lem-mono-0}
T(X)>1, \qquad \mbox{where}\quad T(X)\equiv f_1(X)/f_0(X).  
\eeq
We then use this test for testing $H_2: X_i|Q\overset{iid}{\sim}N(0, \Sigma_2(Q))$ against $H_0$.  
The type I error remains the same as before. Let $P_\ell$ be the probability measure associated with the model $H_\ell$, $\ell=1,2$. If we can show that the type II error does not increase, i.e., 
\beq \label{lem-mono-1}
P_2(T(X)\leq 1) \leq P_1(T(X)\leq 1),
\eeq
then the sum of type I and type II errors does not increase either, i.e., it is no larger than $1-(1/2)\|f_1-f_0\|_1$. On the other hand, by Neyman-Pearson lemma again, with the test \eqref{lem-mono-0} for testing $H_2$ against $H_0$, the sum of type I and type II errors is no smaller than $1-(1/2)\|f_2 - f_0\|_1$. This proves that $\|f_1-f_0\|_1\leq \|f_2 - f_0\|_1$. 

It remains to show \eqref{lem-mono-1}. Without loss of generality, we assume $\Sigma_0(Q)$ is non-degenerate almost surely.  Write $\hat{\Sigma}=(1/n)\sum_{i=1}^n X_iX_i'$ and let $\lambda_1\geq \lambda_2\geq \cdots \geq\lambda_p\geq 0$ be its eigenvalues. We first show that $T(X)$ depends on $X$ only through $\lambda_1,\cdots, \lambda_p$, i.e., 
\beq \label{lem-mono-3}
T(X)=T(\lambda_1,\cdots,\lambda_p). 
\eeq
Let $\hat{\Sigma}=U\Lambda U'$ be the eigen-decomposition of $\hat{\Sigma}$, where $\Lambda=\diag(\lambda_1,\lambda_2,\cdots,\lambda_p)$ and $U$ is an orthogonal matrix. For $\ell=0,1,2$, write $\Omega_\ell(Q)=[\Sigma_\ell(Q)]^{-1}$ and $\tilde{\Omega}_\ell(Q)=U'\Omega_\ell(Q)U$. By our assumption, for any $U$, $\Omega_\ell(Q)$ and $\tilde{\Omega}_\ell(Q)$ have the same distribution. 
It follows that  
\begin{align*}
f_\ell(X) &= \frac{1}{(2\pi)^{np/2}}\int\exp\Big( -\frac{n}{2}\tr(U\Lambda U'\Omega_{\ell}(Q)) +\frac{n}{2}\log(|\Omega_{\ell}(Q)|)  \Big)dF(Q)\cr
&= \frac{1}{(2\pi)^{np/2}} \int\exp\Big(- \frac{n}{2}\tr(\Lambda \tilde{\Omega}_\ell(Q)) +\frac{n}{2}\log(|\tilde{\Omega}_{\ell}(Q)|)  \Big)dF(Q)\cr
&= \frac{1}{(2\pi)^{np/2}} \int\exp\Big(- \frac{n}{2}\tr(\Lambda \Omega_{\ell}(Q)) +\frac{n}{2}\log(|\Omega_{\ell}(Q)|)  \Big)dF(Q), 
\end{align*}
where the second equality is because $|\Omega_\ell(Q)|=|U\tilde{\Omega}_\ell(Q)U'|=|\tilde{\Omega}_\ell(Q)|$. This implies that $f_{\ell}(X)$ depends on $X$ only through $\Lambda$, i.e., $f_{\ell}(X)=f_{\ell}(\lambda_1,\cdots,\lambda_p)$, for $\ell=0,1,2$. 
Combining it with the definition of $T(X)$ gives \eqref{lem-mono-3}. 

From the above expression of $f_\ell(X)$, we see that $\frac{\partial}{\partial\lambda_j}f_\ell$ exists for $\ell=0,1,2$, so $\frac{\partial}{\partial\lambda_j}T$ is well defined. We next show that 
\beq \label{lem-mono-2}
\frac{\partial T}{\partial \lambda_j} \geq 0, \qquad \mbox{for all }1\leq j\leq p. 
\eeq
Fix $1\leq j\leq p$. We note that $\frac{\partial T}{\partial \lambda_j}$ has the same sign as  $f_0\frac{\partial f_1}{\partial \lambda_j}-f_1\frac{\partial f_0}{\partial \lambda_j}$. Define $h(B;\Lambda)= (2\pi)^{-np/2}e^{-\frac{n}{2}\tr(\Lambda B) +\frac{n}{2}\log(|B|)}$ . Then, 
\begin{align*}
& f_\ell(\lambda_1,\cdots,\lambda_p)=E_Q\bigl[ h\big(\Omega_\ell(Q); \Lambda\big)\bigr],\cr
&\frac{\partial f_\ell}{\partial\lambda_j}(\lambda_1,\cdots,\lambda_p)= - \frac{n}{2}E_Q\bigl[ \Omega^{jj}_{\ell}(Q) \cdot h\big(\Omega_\ell(Q);\Lambda\big)\bigr],
\end{align*}
where $\Omega^{jj}_{\ell}(Q)$ is the $j$-th diagonal of $\Omega_\ell(Q)$ and $E_{Q}[\cdot]$ denotes the expectation with respect to the randomness of $Q$ only. Let $Q$ and $\widetilde{Q}$ be two independent copies of $Q$. 
Write $h(B;\Lambda)=h(B)$ for short. 
It is seen that 
\begin{align*}
f_0\frac{\partial f_1}{\partial \lambda_j}
&=  E_Q \bigl[h(\Omega_0(Q))\bigr]\cdot - \frac{n}{2}E_{\tilde{Q}}\bigl[\Omega^{jj}_{1}(\widetilde{Q})\cdot h(\Omega_1(\widetilde{Q}))\bigr]\cr
& = - \frac{n}{2} E_{Q,\tilde{Q}}\bigl[\Omega^{jj}_{1}(\widetilde{Q})\cdot h(\Omega_1(\widetilde{Q}))\cdot h(\Omega_0(Q))\bigr]\cr
& \geq - \frac{n}{2} E_{Q,\tilde{Q}} \bigl[\Omega^{jj}_{0}(Q)\cdot h(\Omega_1(\widetilde{Q}))\cdot h(\Omega_0(Q))\bigr]\cr
&= E_{\tilde{Q}}\bigl [h(\Omega_1(\widetilde{Q}))\bigr]\cdot - \frac{n}{2}  E_Q\bigl[ \Omega^{jj}_{0}(Q)\cdot h(\Omega_0(Q))\bigr]=f_1\frac{\partial f_0}{\partial \lambda_j},
\end{align*}
where the inequality is due to the assumption that $\Sigma_1(\tilde{Q})\succeq \Sigma_0(Q)$ for any realized $Q, \tilde{Q}$. This proves \eqref{lem-mono-2}.

We now show \eqref{lem-mono-1}. 
For $\ell=0,1,2$, let $\lambda_k^{(\ell)}(Q)$ be the $k$-th leading eigenvalue of $\hat{\Sigma}$ under the model $H_{\ell}$. By our assumption, $\Sigma_1(Q)\preceq \Sigma_2(Q)$ for any realized $Q$, which implies $\lambda^{(1)}_k(Q)\leq \lambda_k^{(2)}(Q)$, $1\leq k\leq p$. Combining this with  \eqref{lem-mono-3}-\eqref{lem-mono-2}, 
we have $T(\lambda_1^{(1)}(Q),\cdots,\lambda_p^{(1)}(Q))\leq T(\lambda_1^{(2)}(Q),\cdots,\lambda_p^{(2)}(Q))$. 
It follows that  
\[
P_2(T(X)>1|Q)\geq P_1(T(X)>1|Q), \qquad \mbox{for any realized } Q. 
\]
Taking expectation with respect to $Q$ on both sides gives  \eqref{lem-mono-1}. 
\qed

\subsection{Proof of Theorem~\ref{thm:critical1}} 
By Lemmas~\ref{lem:differ}, \ref{lem:approx} and \ref{lem:tildeR}, 
\begin{align*}
\log &(LR_n) = \log(\widetilde{R}_n) +o_P(1) \cr
&= -\frac{nr^2\delta^2}{4p(1-\delta)^2} -  \frac{nr\delta}{2} \left[\int \psi_n(\lambda) \mu_n^X(d\lambda) - \int \psi_n(\lambda) \mu_{n,p}(d\lambda)\right] +o_P(1), 
\end{align*}
where $\psi_n(\lambda)=\frac{p}{n\delta}\log(1- \frac{\delta n}{p}(\lambda-\frac{1}{1-\delta}))$. With high probability, the support of $\mu_n^X$ and the support of $\mu_{n,p}$ are both strictly contained in a bounded interval $[(1-\sqrt{\gamma})^2-0.1, (1+\sqrt{\gamma})^2+0.1]$. Write $\gamma_n=p/n$. 
Since $\delta=o(1)$, the Taylor expansion yields $\psi_n(\lambda) = \frac{\gamma_n}{\delta}[-\frac{\delta}{\gamma_n}(\lambda-\frac{1}{1-\delta})+O(\delta^2)] = -(\lambda - \frac{1}{1-\delta}) + O(\delta)$. Moreover, since $\int\mu_n^X(d\lambda)=\int\mu_{n,p}(d\lambda)=1$, the term $\frac{1}{1-\delta}$ does not affect the integral. As a result,
\begin{align*}
\log(LR_n) &= - \frac{nr^2\delta^2}{4p(1-\delta)} + \frac{nr\delta}{2}\big[\int \lambda \mu_n^X(d\lambda) - \int \lambda\mu_{n,p}(d\lambda)\big] + o_P(1)\cr
&= - \frac{(r\delta)^2}{4\gamma_n} + \frac{(r\delta)}{2\gamma_n}[S_n - p] + o_P(1), 
\end{align*}
where $S_n$ is the trace statistic, and the second line is because $\int \lambda \mu_n^X(d\lambda)=p^{-1}S_n$ and $\int\lambda\mu_{n,p}(d\lambda)=1$. In other words, 
\[
\log(LR_n) = -\frac{\theta_n^2}{2}  - \theta_n \cdot \left(\frac{S_n - p}{\sqrt{2\gamma_n}}\right) +o_P(1), \qquad \mbox{where}\;\; \theta_n\equiv \frac{r\delta}{\sqrt{2\gamma_n}}\goto \theta.  
\]
Note that $\frac{S_n-p}{\sqrt{2\gamma_n}} =\frac{(nS_n-np)}{\sqrt{2np}}$. By elementary statistics, this term converges to $N(0,1)$ under the null, and converges to $N(\frac{nr\delta}{\sqrt{2np}},1)=N(\theta_n, 1)$ under the alternative. The claim then follows. \qed

\subsection{Proof of Lemma~\ref{lem:differ}}  \label{subsec:proof-differ}
Let $Y\in \mathbb{R}^{p,r}$ be a matrix such that $Y(i,j)\overset{iid}{\sim}N(0,1)$. It is known that when $Q$ is uniformly distributed over the Stiefel manifold $\mathbb{S}(p,r)$, $QQ'$ has the same distribution as $Y(Y'Y)^{-1}Y'$. So, the RWS model \eqref{model-alt} is equivalent to 
\beq \label{thm-2model-0}
H_1^{(n)}: \quad \Sigma = [I_p - \delta Y(Y'Y)Y']^{-1}. 
\eeq
Denote by $\eta_1\geq \eta_2\geq \cdots\geq\eta_r>0$ the nonzero eigenvalues of $(1/p)YY'$. Let $B_n=\{Y: \max_{1\leq k\leq r}|\eta_k-1|\leq 2\sqrt{r/p}\}$.

We first show that $Y\in B_n^c$ has a negligible effect on the likelihood ratio. Denote by $\overline{\Sigma}(Y)$ a general covariance matrix model indexed by $Y$. Introduce
\beq \label{thm-2model-tmp}
\overline{\Sigma}^*(Y)=\overline{\Sigma}(Y)\cdot 1\{Y\in B_n\} + I_p\cdot 1\{Y\in B_n^c\}. 
\eeq
Consider the model $X_i|Y\overset{iid}{\sim}N(0, \overline{\Sigma}(Y))$ and  let $\overline{LR}_n$ be the corresponding likelihood ratio with respect to $H_0^{(n)}: X_i\overset{iid}{\sim} N(0, I_p)$. Define $\overline{LR}_n^*$ similarly by replacing $\overline{\Sigma}(Y)$ with $\overline{\Sigma}^*(Y)$. Below, we show that 
\beq  \label{smallevent}
E_0|\overline{LR}_n - \overline{LR}_n^*| \goto 0. 
\eeq
To see this, let $f_0(X)$ be the joint density of $X$ under $H_0^{(n)}$, $f(Y)$ be the joint density of $Y$ and $f_1(X|Y)$ be the conditional density of $X$ under $X_i|Y\overset{iid}{\sim}N(0, \overline{\Sigma}(Y))$. Then, $\overline{LR}_n=\int \frac{f_1(X|Y)f(Y)}{f_0(X)}dY$ and $\overline{LR}_n^* = \int_{B_n} \frac{f_1(X|Y)f(Y)}{f_0(X)}dY+\int_{B_n^c}f(Y)dY$. 
So $E_0|LR_n-LR_n^*|\leq \int f_0(X)\int_{B_n^c}\big[\frac{f_1(X|Y)f(Y)}{f_0(X)}+ f(Y)\big] dYdX =\int_{B_n^c} [\int f_1(X|Y)dX]f(Y) dY + \int_{B_n^c}[\int f_0(X)dX] f(Y)dY = 2\int_{B_n^c}f(Y)dY= 2\cdot P(Y\in B_n^c)$. By \cite[Corollary 5.35]{Vershynin}, $P(Y\in B_n^c)\leq 2e^{-r}$. So \eqref{smallevent} follows. 

Regarding \eqref{smallevent}, for all the covariance models considered here, switching between $Y$ and $Y\cdot 1\{Y\in B_n\}$ only affects the $L^1$-distance by an $o(1)$ term, and we omit such a difference for notation simplicity.

Introduce $\delta^{\pm}=\delta/(1\mp 2\sqrt{p/n})$ and let $\widetilde{LR}_n^{\pm}$ be the likelihood ratio associated with the Gaussian proxy model with $\delta=\delta^{\pm}$. For $Y\in B_n$, 
\beq \label{thm-2model-matrices}
\delta^- (p^{-1} YY')\preceq \delta [Y(Y'Y)^{-1}Y'] \preceq \delta^+ (p^{-1} YY'). 
\eeq
We apply Lemma~\ref{lem:monotone2} with $\Sigma_0(Y)=I_p$, $\Sigma_1(Y)=[I_p-\delta Y(Y'Y)^{-1}Y']^{-1}$ and $\Sigma_2(Y)=[I_p - \delta^+ p^{-1} YY']^{-1}$.\footnote{When applying Lemma~\ref{lem:monotone2}, we consider $Y\cdot 1\{Y\in B_n\}$ in all the covariance models. So, we need to show that for any nonradom orthogonal matrix $U\in\mathbb{R}^{p,p}$, $UY\cdot1\{UY\in B_n\}$ has the same distribution as $Y\cdot 1\{Y\in B_n\}$. We note that $Y\overset{(d)}{=}P_1\Gamma P_2$, where $P_1$ and $P_2$ are uniformly distributed over $\mathbb{S}(p,r)$ and $\mathbb{S}(r,r)$, respectively, and $\Gamma\in\mathbb{R}^{r,r}$ is a random diagonal matrix, and $(P_1, \Gamma, P_2)$ are independent. So $UY\cdot 1\{UY\in B_n\}\overset{(d)}{=}UP_1\Gamma P_2\cdot 1\{\|\Gamma-I_r\|\leq 2\sqrt{r/p}\}$. The claim follows by noting that $UP_1$ has the same distribution as $P_1$.} It yields that
\[
E_0|1-\widetilde{LR}_n^-|\leq E_0|1-LR_n|\leq E_0|1-\widetilde{LR}_n^+|. 
\]
Similarly, since $\delta^-(p^{-1}YY')\preceq \delta (p^{-1}YY')\preceq \delta^+ (p^{-1}YY')$, 
\[
E_0|1-\widetilde{LR}_n^-|\leq E_0|1-\widetilde{LR}_n|\leq E_0|1-\widetilde{LR}_n^+|. 
\]
As a result,
\[
| E_0|1-\widetilde{LR}_n| - E_0|1-LR_n| | \leq E_0|1-\widetilde{LR}^+_n| - E_0|1-\widetilde{LR}^-_n|\leq E_0| \widetilde{LR}^-_n -  \widetilde{LR}^+_n|. 
\]
To show the claim, it suffices to show that
\beq  \label{thm-2model-1}
E_0| \widetilde{LR}^-_n -  \widetilde{LR}^+_n| \goto 0. 
\eeq

We now show \eqref{thm-2model-1}. Denote by $\mu_n^X$ the empirical spectral measure of $\hat{\Sigma}$, and by $\mu_{n,p}$ the Marchenko-Pastur (MP) measure with parameter $\gamma=p/n$ (see Section~\ref{sec:RMT}). Let $F_n^X$ and $F_{n,p}$ be the corresponding CDF's. Introduce $E_n=\{X\in\mathbb{R}^{n,p}: \sup_{t\in\mathbb{R}}|F_n^X(t)-F_{n,p}(t)|\leq n^{-1}\log^6(n)\log^6(\log(n))\}$. By Theorem~\ref{thm:KS}, $P_0(X\in E_n^c)=O(e^{-c\log(n)\log(\log(n))})$, as $n\goto\infty$. Note that
\[
E_0\bigl(| \widetilde{LR}^-_n -  \widetilde{LR}^+_n|\cdot I_{E_n^c}\bigr) \leq E_0 (\widetilde{LR}^-_n\cdot I_{E_n^c})
+E_0 (\widetilde{LR}^+_n\cdot I_{E_n^c}). 
\]
Since $E_0(\widetilde{LR}^{\pm})=1$ is finite, $P_0(E_n^c)=o(1)$ implies $E_0 (\widetilde{LR}^-_n\cdot I_{E_n^c})=o(1)$. Therefore, to show \eqref{thm-2model-1}, it suffices to show that
\beq  \label{thm-2model-1(a)}
E_0\bigl(| \widetilde{LR}^-_n -  \widetilde{LR}^+_n|\cdot I_{E_n}\bigr) \goto 0. 
\eeq

Below, we show \eqref{thm-2model-1(a)}. The assumption on $(\alpha,\beta)$ ensures that $e_n \leq n^{-\epsilon}$ for some $\epsilon>0$, where
\[
e_n \equiv \delta^2 r^{5/2}n^{-1/2} +  \delta r^{3/2}n^{-1/2} + \delta^2 r.  \footnote{The first two terms are from the calculation of $(I)$ and $(II)$ below, and the last term is due to $E_0|\widetilde{LR}_n-\widetilde{R}_n|\leq C(\delta^2 r + \delta^3 r^{5/2} n^{-1/2})$, as seen in the proof of Lemma~\ref{lem:approx}.}
\]
In Section~\ref{subsec:proxy}, we define $\widetilde{R}_n$, a proxy of  $\widetilde{LR}_n$. Let $\widetilde{R}^{\pm}_n$ be the proxy associated with $\delta=\delta^{\pm}$. By Lemma~\ref{lem:approx}, $E_0|\widetilde{LR}^{\pm}_n-\widetilde{R}^{\pm}_n|=o(1)$. So it suffices to show
\beq  \label{thm-2model-2}
E_0\bigl(| \widetilde{R}^-_n -  \widetilde{R}^+_n|\cdot I_{E_n}\bigr) \goto 0. 
\eeq
Let $g_n^{\pm}(\lambda)=\frac{p}{n}\log(1- \frac{\delta^{\pm} n}{p}(\lambda-\frac{1}{1-\delta^{\pm}}))$. By Lemma~\ref{lem:tildeR}, 
\[
\log(\widetilde{R}_n^{\pm} ) = -\frac{nr^2(\delta^{\pm})^2}{4p(1-\delta^{\pm})^2} -  \frac{nr}{2} \left[\int g_n^{\pm}(\lambda) \mu_n^X(d\lambda) - \int g_n^{\pm}(\lambda) \mu_{n,p}(d\lambda)\right]. 
\]
Let $\Delta g_n = g_n^{+} - g_n^{-}$. It follows that
\begin{align*}
&\log\bigl(\widetilde{R}_n^{-}/\widetilde{R}_n^+\bigr) = \frac{nr}{2} \left[\int \Delta g_n(\lambda) \mu_n^X(d\lambda) - \int \Delta g_n(\lambda) \mu_{n,p}(d\lambda)\right]  \cr
& + \frac{nr^2(\delta^{+})^2}{4p(1-\delta^{+})^2} -\frac{nr^2(\delta^{-})^2}{4p(1-\delta^{-})^2}\equiv (I) + (II). 
\end{align*}
On the event $E_n$, 
\[
|(I)|\leq \frac{nr}{2}\cdot \sup_{\lambda}|(\Delta g_n)'(\lambda)|\cdot \sup_{\lambda\in\mathbb{R}}|F_n^X(\lambda)-F_{n,p}(\lambda)|\leq L_n r\sup_{\lambda}|(\Delta g_n)'(\lambda)|,
\]
where $L_n=\log^6(n)\log^6(\log(n))/2$. For $X\in E_n$, the supports of $\mu_{n,p}$ and $\mu_n^X$ are both strict subsets of  $A=[\frac{1}{2}(1-\sqrt{\gamma})^2, 2(1+\sqrt{\gamma})^2]$. So it suffices to consider $\sup_{\lambda\in A}|(\Delta g_n)'(\lambda)|$. By direct calculations, $(\Delta g_p)'(\lambda)=-(\delta^+-\delta^{-})\cdot [1- \tfrac{n\delta^{+}\delta^{-}}{p(1-\delta^{+})(1-\delta^{-})}][1-\tfrac{\delta^{+}n}{p}(\lambda-\tfrac{1}{1-\delta^{+}})]^{-1}[1-\tfrac{\delta^{-}n}{p}(\lambda-\tfrac{1}{1-\delta^{-}})]^{-1}=- (\delta^{+}-\delta^{-})[1+o(1)]$, for $\lambda\in A$.  It follows that
\[
\sup_{\lambda\in A}|(\Delta g_n)'(\lambda)|\leq C(\delta^{+}-\delta^{-})\leq C\delta\sqrt{r/p}. 
\]
Combining the above gives
\beq \label{thm-2model-3}
|(I)|\leq CL_n\cdot \delta r^{3/2}n^{-1/2}. 
\eeq
By direct calculation, 
$(II) = \frac{nr^2}{4p}\frac{(\delta^+-\delta^{-})(\delta^{+}+\delta^{-}-2\delta^+\delta^-)}{(1-\delta^{+})^2(1-\delta^{-})^2}$. It follows that 
\beq \label{thm-2model-4}
|(II)| \leq Cr^2 (\delta^+-\delta^{-})\delta \leq C\delta^2 r^{5/2}n^{-1/2}. 
\eeq
Both the right hands of  \eqref{thm-2model-3}-\eqref{thm-2model-4} are $O(e_n)$. As a result, for $X\in E_n$, 
\[
|\widetilde{R}_n^{-}/\widetilde{R}_n^{+} -1|\leq \max\{|e^{Ce_n}-1|, |e^{-Ce_n}-1|\}\leq Ce_n.  
\]
So $E_0(| \widetilde{LR}^-_n -  \widetilde{LR}^+_n|\cdot I_{E_n})\leq Ce_nE_0(\widetilde{R}_n^+)$. Since $E_0|\widetilde{R}_n^+-\widetilde{LR}_n^+|=o(1)$ and $E_0(\widetilde{LR}_n^+)=1$, $E_0(\widetilde{R}_n^+)\leq 2$. Combining these gives \eqref{thm-2model-2}. 
\qed

{\bf Remark}. ({\it Conjecture on $E_0|LR_n-\widetilde{LR}_n|$}).  We conjecture that a stronger result holds:
\[
E_0|LR_n-\widetilde{LR}_n|\goto 0. 
\]
To see this, let $\Sigma(Y)$ be the covariance matrix \eqref{thm-2model-0} of RWS, and let $\widetilde{\Sigma}(Y)$, $\widetilde{\Sigma}^{\pm}(Y)$ be the covariance matrices of proxy models associated with $\delta$, $\delta^{\pm}$, respectively. By \eqref{thm-2model-matrices}, $\widetilde{\Sigma}^{-}(Y)\preceq \Sigma(Y)\preceq \widetilde{\Sigma}^{+}(Y)$, so the RWS is ``sandwiched" by two proxy models. We conjecture that 
\beq \label{thm-2model-guess}
E_0|LR_n-\widetilde{LR}_n^-|\leq E_0|\widetilde{LR}_n^+ - \widetilde{LR}_n^-|. 
\eeq
Similarly, since $\widetilde{\Sigma}^{-}(Y)\preceq \widetilde{\Sigma}(Y)\preceq \widetilde{\Sigma}^{+}(Y)$, we conjecture $E_0|LR_n-\widetilde{LR}_n^-|\leq E_0|\widetilde{LR}_n^+ - \widetilde{LR}_n^-|$. By the triangular inequality, $E_0|LR_n-\widetilde{LR}_n|\leq E_0|LR_n-\widetilde{LR}_n^-| + E_0|\widetilde{LR}_n-\widetilde{LR}_n^-|\leq 2E_0|\widetilde{LR}_n^+ - \widetilde{LR}_n^-|$. It then follows from \eqref{thm-2model-1} that $E_0|LR_n-\widetilde{LR}_n|\goto 0$. 

To rigorously prove \eqref{thm-2model-guess}, we need a slightly stronger result than Lemma~\ref{lem:monotone2} (monotonicity of $L^1$-distance). In Lemma~\ref{lem:monotone2}, one assumption is that $\Sigma_0(Q)\preceq \Sigma_1(\widetilde{Q})$ for any realizations $Q,\widetilde{Q}$. If this assumption is relaxed to requiring $\Sigma_0(Q)\preceq \Sigma_1(Q)$ for any $Q$, then we can apply the result to $(\Sigma_0, \Sigma_1, \Sigma_2)=(\widetilde{\Sigma}^-(Y),\Sigma(Y),\widetilde{\Sigma}^+(Y))$ and obtain \eqref{thm-2model-guess}.

\section{Discussions}  \label{sec:discuss}


Testing of sphericity, i.e., whether the population covariance matrix $\Sigma$ is the identity,
has been well-studied. While many tests were proposed for high-dimensional settings \cite{bai2009corrections, ChenSX, Johnstone2001, ledoit2002some}, there are relatively few results about the phase transition: \cite{cai2013optimal} studied the testing limit for a class of alternatives that are characterized by the Frobenius norm of $\Sigma-I_p$, and \cite{BR13} investigated a class of alternatives which are rank-$1$ sparse perturbations of the null. We are interested in a different class of alternatives, the Rare and Weak Spike model. When the number of spikes is finite, \cite{Hallin2013,Hallin2014} studied the asymptotic power envelopes of eigenvalue-based tests, and their results implicitly gave the phase transition of this testing problem. Compared with \cite{Hallin2013,Hallin2014}, our results cover broader situations where the number of spikes may tend to infinity as $n,p\goto\infty$. We also extend our results to more complicated settings where $\Sigma$ is not necessarily the identity in the null, by exploring the idea of ``comparison of experiments." 

Various tests were proposed for testing $\Sigma=I_p$; some are based on extreme eigenvalues \cite{Johnstone2001,Johnstone2015}, and some are based on bulk eigenvalues \cite{bai2009corrections,dobriban2016sharp, ledoit2002some, Hallin2013}. It has been recognized that the extreme-eigenvalue tests are more powerful for testing Rare and Strong spikes, while for Rare and Weak spikes, it is better to use more than just a few top eigenvalues. In the spirit of ``letting the data decide which eigenvalues to use", we proposed two new tests, the Higher Criticism (HC) test, and the Cumulative Sum (CuSum) test. Both tests can be viewed as extensions of the Higher Criticism test for testing sparse normal means \cite{DJ04}. In contrast to bulk-eigenvalue tests \cite{dobriban2016sharp,Hallin2013}, our proposed tests are purely adaptive. Moreover, we find that the CuSum test is both theoretically optimal and has impressive numerical performance. 

Technically, we introduce a Gaussian proxy model for studying the RWS model. We show that the two models are close enough for a wide range of parameters. At the same time, the likelihood of the proxy model is easier to analyze, where we can use the available results on Wishart matrices instead of analyzing any spherical integrals. Such an approach can be useful for other problems related to the RWS model.

In Random Matrix Theory, it is an interesting topic to study the empirical eigenvalues of a spiked covariance model. In the seminal work \cite{BAP}, they derived the limiting joint distribution of top $k$ eigenvalues of a complex Gaussian Wishart matrix, with the population covariance a rank-$r$ perturbation of the identity, where $r$ is finite and $k\leq r$. 
Later, \cite{Baik, Paul2007} proved similar results for real covariance matrices. However, there is little understanding on the case that $r$ tends to infinity as $n,p$ grows. Moreover, existing works mainly focus on a few top eigenvalues. The bulk eigenvalues of a spiked model are also of great interest in many applications, but such results are only found for the identify covariance case \cite{Sandrine2013,PillaiYin, Su}.
Our results fill in the gap by both allowing $r\goto\infty$ (in a algebraic rate) and concerning all eigenvalues.  

A somewhat relevant setting to the spiked covariance model is the spiked Wigner model, where we consider an $n\times n$ Wigner matrix plus a non-random, rank-$r$ diagonal matrix. 
Such a setting with $r\goto\infty$ was studied in \cite{peche2006largest}, and they derived the limiting behavior of top eigenvalues. Besides that the two settings are very different, \cite{peche2006largest} assumes the magnitude of entries in the diagonal perturbation is independent of $n$, so it is not about rare and ``weak" spikes. 

Our work also brings a new angle of using statistical theory to shed lights on random matrix theory. We recognize that, for parameters where all tests are asymptotically powerless, any bounds of empirical eigenvalues that hold in the null should also hold in the alternative. Hence, as long as we know the statistical limits, we can take advantage of results in the identity covariance case for studying the spiked covariance models. Alternatively, if we wish to derive such results directly from random matrix theory, we have to pay non-trivial efforts, noting that how to extend \cite{BAP,Baik, Paul2007} to the setting of rare and weak spikes is largely unclear.


It is possible to extend our results to nonGaussian settings: $X=Y\Sigma^{1/2}$, where $Y$ is an $n\times p$ random matrix with $iid$ entries of mean $0$ and variance $1$. To test $\Sigma=I_p$, we can still apply the four tests considered in this paper, and their performance will be similar to the Gaussian case, due to universality of eigenvalue statistics (e.g., \cite{PillaiYin}). The lower bound arguments may be slightly different, depending on the distribution of the entries of $Y$. However, the key technique of this paper continues to work:  we can approximate the RWS model (on $\Sigma$) by a Gaussian proxy model, to overcome the hurdle of analyzing likelihood ratios. We leave this for future work. 

Another extension of our work is to consider the same testing problem for population {\it correlation} matrices. For this setting, it is natural to apply the tests to eigenvalues of the sample correlation matrix. The existing results about the null behavior of these eigenvalues \cite{bao2012tracy,elkaroui2009} will be helpful. We also leave this for future study. 


\section{Appendix} 
\label{sec:appen}

\subsection{Proof of Theorem~\ref{thm:deviation}}
Define the function $h:\mathbb{R}^{n,p}\to \mathbb{R}$ by $h(Z)=p\mu_n^Z(f)$, i.e., we view $p\mu_n^Z(f)=\sum_{j=1}^p f(\lambda_j)$ as a function of the entries of $Z$. 
For a differentiable $f$ such that $\|f'\|_\infty<\infty$, Delyon \cite{Delyon2010} shows that $\nabla h$ exists. According to the last equation on \cite[Page 554]{Delyon2010}, for each $1\leq i\leq n$, $\sum_{1\leq j\leq p} \big(\frac{\partial h(Z)}{\partial Z_{ij}}\big)^2\leq 4n^{-2}\|f'\|_\infty^2 \sum_{1\leq j\leq p}Z_{ij}^2$. 
It follows that 
\beq \label{thm-deviation-1}
\|\nabla h(Z)\|^2\leq 4n^{-2}\|f'\|_\infty^2\|Z\|_F^2.
\eeq 
Introduce
\[
\tilde{h}(Z)=\frac{\sqrt{n\log(2)}}{2\sqrt{2p}\|f'\|_\infty}h(Z). 
\]
Since $\|Z\|_F^2$ has a $\chi^2_{np}(0)$ distribution, $E[e^{t\|Z\|_F^2}]=(1-2t)^{-np/2}$ for all $t\geq 0$. From \eqref{thm-deviation-1}, $\|\nabla \tilde{h}(Z)\|^2\leq \frac{\log(2)}{2np}\|Z\|_F^2$. It follows that 
\beq \label{thm-deviation-2}
E[e^{\|\nabla \tilde{h}(Z)\|^2}]\leq \bigl[1-(2np)^{-1}\log(2)\bigr]^{-np/2}\leq 2. 
\eeq
We now apply \cite[Theorem 6.1]{Bobkov2015}, a concentration inequality for non-Lipschitz functions, to $\tilde{h}$. The probability measure $\mu$ there corresponds to the standard Gaussian measure of dimension $np$ in our settings, so the subGaussian constant $\sigma(\mu)=1$. It yields that for $L_0$ such that $P(|\nabla\tilde{h}(Z)|>L_0)\leq 1/2$ and all $t>0$, 
\[
P\big(  |\tilde{h}(Z)-E[\tilde{h}(Z)]|>t  \big)\leq 2\inf_{y\leq L_0} \big\{  e^{-\frac{t^2}{cy^2}} +   P(  \|\nabla\tilde{h}(Z)\|>y ) \big\},
\]
where $c>0$ is an absolute constant. By \eqref{thm-deviation-2} and the Chebyshev inequality, $P(  \|\nabla\tilde{h}(Z)\|>L )\leq 2e^{-L^2}$. Therefore, for any $y\geq L_0\equiv \sqrt{2\log(2)}$, 
\beq \label{thm-deviation-3}
P\big(  |\tilde{h}(Z)-E[\tilde{h}(Z)]|>t  \big) \leq 2e^{-\frac{t^2}{cy^2}} + 4e^{-y^2}. 
\eeq
The claim follows by noticing that $\mu_n^Z(f) = \frac{2\sqrt{2}\|f'\|_\infty}{\sqrt{np\log(2)}}\tilde{h}(Z)$. \qed

\subsection{Proof of Lemma~\ref{lem:eigen-moment}}
Fix $m\geq 1$. Write for short $H=\sum_{k=1}^p (\lambda_j-1)^m$ and $e_n=\log(n)[(p\vee\log(n))/n]^{m/2}$. We aim to show that for some constant $C>0$,   
\[
P(|H-E(H)|> Ce_n) = o(n^{-1}). 
\]
Note that $p\leq n$. By \cite[Corollary 5.35]{Vershynin}, for every $t\geq 0$,
\beq \label{Vershynin5.35}
P\left(   \max_{1\leq j\leq p}|\lambda_j-1| > \sqrt{p/n} +t/\sqrt{n} \right)\leq 2e^{-t^2/2}. 
\eeq
Let $a_n=\sqrt{p/n}+2\sqrt{\log(n)/n}$ and $B_n$ be the event $\max_{1\leq j\leq p}|\lambda_j-1|\leq a_n$. The above implies that $P(B_n^c)=o(n^{-1})$. So it suffices to show that 
\beq \label{lem-mmt-1}
P( |H-E(H)|>Ce_n,\; B_n ) = o(n^{-1}). 
\eeq
Define a function $f_n(x)=h_n(x-1)$, where 
\[
h_n(x) = 
\begin{cases}
x^m, & \mbox{if } 0\leq x\leq a_n,\\
(2a_n)^m, & \mbox{if } x\geq 2a_n,\\
a_n^m + m\int_{0}^{x-a_n} (2a_n-x)^{m-1}, & \mbox{if }a_n<x<2a_n, \\
(-1)^m h_n(-x), & \mbox{if }x<0. 
\end{cases} 
\]
Introduce $H^* = \sum_{j=1}^p f_n(\lambda_j)$. Then, $H=H^*$ over the event $B_n$. Therefore, to show \eqref{lem-mmt-1}, it suffices to show that for some constants $C_1, C_2>0$, 
\beq  \label{lem-mmt-2}
P(|H^* - E(H^*)|> C_1e_n) = o(n^{-1}), 
\eeq
and
\beq \label{lem-mmt-3}
|E(H^*) - E(H)| < C_2 e_n. 
\eeq

First, consider \eqref{lem-mmt-2}. We write $H^*=p\mu_n^Z(f_n)$, and note that $f_n$ is continuously differentiable and $\|f'_n\|_\infty=ma_n^{m-1}$. By Theorem~\ref{thm:deviation} where we set $y=\sqrt{2}\log(n)$ and $t=2\sqrt{c_0}\log(n)$, with probability at least $1-o(n^{-1})$,
\[
|H^* - E(H^*)|\leq \frac{2\sqrt{c_0p}\log(n)}{\sqrt{n}} \|f'_n\|_\infty\leq Cma_n^{m-1}\sqrt{p/n}\log(n). 
\]
The right hand side is no larger than $Ce_n$. This proves \eqref{lem-mmt-2}. 

Next, consider \eqref{lem-mmt-3}. Let $\xi=\max_{1\leq j\leq p}|\lambda_j-1|$. We have
\beq \label{lem-mmt-4}
|E(H)-E(H^*)| \leq \sum_{j=1}^p E(|\lambda_j-1|^m \cdot I_{B_n^c})\leq pE\bigl(\xi^m \cdot I_{\{\xi>a_n\}}\bigr). 
\eeq
We use the simple fact that for any continuous random variable $W$ and real number $a$, $E(W\cdot I_{\{W>a\}})=aP(W>a)+ \int_{a}^{\infty}P(W>t)dt$. It follows that 
\begin{align*}
&E\bigl(\xi^m \cdot I_{\{\xi>a_n\}}\bigr) = a_n^m P( \xi>a_n ) + \int_{a_n^m}^\infty P(\xi^m >x)dx\cr
&= o(n^{-1}a_n^m ) + \frac{m}{\sqrt{n}} \int_{2\sqrt{\log(n)}}^\infty (\sqrt{p/n}+t/\sqrt{n})^{m-1} P\bigl(\xi>\sqrt{p/n}+t/\sqrt{n}\bigr)dt \cr
&= o(n^{-1}a_n^m) +  \frac{m}{(\sqrt{n})^m} \int_{2\sqrt{\log(n)}}^\infty 2^{m-1}[(\sqrt{p})^{m-1}+t^{m-1}]  e^{-t^2/2}dt\cr
&\leq o(n^{-1}a_n^m) + \frac{2^{m-1}m}{(\sqrt{n})^m}\cdot \big[ (\sqrt{p})^{m-1}+ (2\sqrt{\log(p)})^{m-1}\big]\cdot Ce^{-2\log(n)}.  
\end{align*}
Here the second line is due to a variable change $t= \sqrt{n}(x^{1/m}-\sqrt{p/n})$, the third line follows from \eqref{Vershynin5.35} and that $(a+b)^{k}\leq 2^k (a^k + b^k)$ for $a,b>0$ and $k$ being a positive integer, and the last line is from elementary calculation. It is easy to see that the right hand side is no larger than $Cp^{-1}e_n$. Combining it with \eqref{lem-mmt-4} gives \eqref{lem-mmt-3}. 
\qed

\subsection{Proof of Lemma~\ref{lem:MPintegral}}
Write $L(t;\gamma)=\int \log(t-\lambda)dF^{mp}_{\gamma}(\lambda)$ for any $\gamma\geq 1$ and $t>(1+\sqrt{\gamma})^2$. It is known that (e.g., equation (3.6) of \cite{Hallin2013}) the Hilbert transformation of $F^{mp}_\gamma$ is 
\[
H(t;\gamma)\equiv \int (t- \lambda)^{-1} dF^{mp}_\gamma (\lambda) = \frac{t + \gamma - 1 - \sqrt{(t - \gamma - 1)^2 - 4 \gamma}}{2 \gamma t},
\]
for $t>(1+\sqrt{\gamma})^2$; the sign of the square root is the same as that of $(t - \gamma - 1)$. 
Since $\frac{\partial}{\partial t}\log(t-\lambda)=(t-\lambda)^{-1}$, we have
\beq \label{lem-MPint-2}
L(t;\gamma) = \int H(t;\gamma)dt + C_0, 
\eeq
where $\int H(t; \gamma)dt$ is any indefinite integral of $H(t;\gamma)$. We take the following choice (it can be verified by checking the derivative):
\begin{align} \label{lem-MPint-3}
\int H(t;\gamma)&dt = \frac{1}{2\gamma}\Big\{ (A_t- B_t) -(\gamma +1)\log(A_t- B_t)\cr
& + (\gamma-1) \log\big((\gamma+1)A_t - (\gamma-1)B_t+4\gamma\big)\Big\}, 
\end{align}
where $A_t = t-\gamma -1$ and  $B_t = \sqrt{(t-\gamma-1)^2-4\gamma}$ are the same as those in Lemma~\ref{lem:MPintegral}. To get the constant $C_0$ in \eqref{lem-MPint-2}, we note that $\lim_{t\to\infty} [L(t;\gamma)-\log(t)]=\lim_{t\to\infty}\int [\log(t-\lambda)-\log(t)]dF^{mp}_{\gamma}(\lambda) = 0$. As a result, 
\[
C_0 = -\lim_{t\goto\infty}[\int H(t;\gamma)dt - \log(t)]. 
\]
As $t\goto\infty$, we have $A_t/t\goto 1$ and $B_t/t\goto 1$. Also, $A_t-B_t=4\gamma/(A_t+B_t)\goto 0$ and $\log(A_t-B_t)+\log(2t)\goto\log(4\gamma)$. Moreover, $(2t)^{-1}[(\gamma+1)A_t - (\gamma-1)B_t+4\gamma]=(2t)^{-1}[A_t+B_t +\gamma(A_t-B_t)+4\gamma]\goto 1$; it follows that $\log((\gamma+1)A_t- (\gamma-1)B_t+4\gamma)-\log(2t)\goto 0$. Plugging these results into \eqref{lem-MPint-3}, we find that
\[
\lim_{t\goto\infty}\Big[\int H(t;\gamma)dt - \frac{\gamma+1}{2\gamma}\log(2t) - \frac{\gamma-1}{2\gamma}\log(2t)\Big] \goto - \frac{\gamma+1}{2\gamma}\log(4\gamma). 
\]
It follows that
\beq \label{lem-MPint-4}
C_0 = - \log(2) + \frac{\gamma+1}{2\gamma}\log(4\gamma). 
\eeq
Combining \eqref{lem-MPint-2}-\eqref{lem-MPint-4} gives the explicit expression of  $L(t;\gamma)$. 
\qed

\subsection{Proof of Lemma~\ref{lem:approx}}
Let $e_n = \delta^3r^{5/2}n^{-1/2} + \delta^2 r$. It suffices to show that when $e_n=o(1)$,  $E_0|\widetilde{LR}_n - \widetilde{R}_n|\leq Ce_n$. 

Let $h(X,Y)$ be the integrand in \eqref{DefinetildeLR}. Define the events $A_n=\{Y: \delta\eta_1<1/2\}$ and $D_n=\{X: p^{-1}\delta n(\lambda_1-\frac{1}{1-\delta})<1\}$. 
Introduce
\[
\Delta_n(Y)=\Delta_n(Y,\delta,r) \equiv g^*_n(Y,\delta, r) - g_n(Y,\delta,r). 
\] 
By \eqref{DefinetildeLR}-\eqref{DefinetildeR}, $\widetilde{LR}_n=\int_{A_n} h(X,Y)dF(Y)+\int_{A_n^c}dF(Y)$ and $\widetilde{R}_n=I_{D_n^c}(X) + I_{D_n}(X)\cdot \int h(X,Y)e^{\frac{n}{2}\Delta_n(Y)}dF(Y)$. 
Let $F_0(X)$ be the joint CDF of $X$ under $H_0^{(n)}$. Then,
\begin{align} \label{approx-1}
E_0 |\widetilde{LR}_n - \widetilde{R}_n| &\leq \int_{A_n\times D_n} h(X,Y)|1-e^{\frac{n}{2}\Delta_n(Y)}|dF(Y)dF_0(X)\cr
&+ \int_{A_n^c\times D_n} |1-h(X,Y)e^{\frac{n}{2}\Delta_n(Y)}|dF(Y)dF_0(X)\cr
&+  \int_{A_n\times D_n^c} |h(X,Y)-1|dF(Y)dF_0(X). 
\end{align}
By definition of the Gaussian model, conditioning on any $Y\in A_n$, $h(X,Y)dF_0(X)$ is the joint density of $X$ under $\widetilde{H}_1^{(n)}$. Consequently,  
\beq  \label{approx-2}
\int h(X,Y)dF_0(X) =1, \qquad \mbox{for}\;\; Y\in A_n. 
\eeq
We plug \eqref{approx-2} into  \eqref{approx-1}, and note that $|1-h(X,Y)|\leq 1+h(X,Y)$ and $|1-h(X,Y)e^{\frac{n}{2}\Delta_n(Y,\delta,r)}|\leq 1+h(X,Y)e^{\frac{n}{2}\Delta_n(Y,\delta,r)}$. It follows that 
\begin{align} \label{approx-3}
&E_0|\widetilde{LR}_n - \widetilde{R}_n| \leq \int_{A_n} |1-e^{\frac{n}{2}\Delta_n(Y)}|dF(Y)\cr
&+ P_0(X\in D_n, Y\in A_n^c) + \int_{A_n^c\times D_n} h(X,Y)e^{\frac{n}{2}\Delta_n(Y)}dF_0(X)dF(Y)\cr
&+ P_0(X\in D_n^c, Y\in A_n) +  \widetilde{P}_1(X\in D_n^c, Y\in A_n),
\end{align}
where $P_0$ and $\widetilde{P}_1$ are the probability measures under $H_0^{(n)}$ and $\widetilde{H}_1^{(n)}$, respectively. By \eqref{Vershynin5.35}, $P_0(X\in D_n^c) = P_0(\lambda_1-1>\frac{\delta}{1-\delta}+\frac{p}{\delta n})\leq 2e^{-\frac{n}{2}(\frac{\delta}{1-\delta}+\frac{p}{\delta n}-\sqrt{p/n})^2}\leq 2e^{-p/2}$, where the last inequality follows from $\frac{\delta}{1-\delta}+\frac{p}{\delta n}\geq 2\sqrt{p/n}$ by Cauchy-Schwartz inequality. 
Again, by \eqref{Vershynin5.35}, $P_0(Y\in A_n^c)\leq 2e^{-\frac{p}{2}[(2\delta)^{-1}-1-\sqrt{r/p}]^2}\leq 2e^{-cn}$. 
Combining these with \eqref{approx-3}, to show the claim, it suffices to show that 
\begin{align}  \label{J1+J2}
J_1 &\equiv  \int_{A_n} |1-e^{\frac{n}{2}\Delta_n(Y)}|dF(Y) \leq Ce_n \cr
J_2 &\equiv \int_{A_n^c\times D_n} h(X,Y)e^{\frac{n}{2}\Delta_n(Y)}dF_0(X)dF(Y) = o(e_n)\cr
J_3 &\equiv \widetilde{P}_1(X\in D_n^c, Y\in A_n) =  o(e_n). 
\end{align}

Consider $J_2$. By \eqref{DefinetildeR}, $h(X,Y)e^{\frac{n}{2}\Delta_n(Y)}= (1-\delta)^{\frac{nr}{2}}b_n(\delta,r)\exp\{\frac{\delta n}{2p}\tr((\hat{\Sigma}-\frac{1}{1-\delta})YY')$. We note that $(1-\delta)^{\frac{nr}{2}}\leq 1$ and   $b_n(\delta,r)\leq e^{\frac{nr\delta}{2p(1-\delta)}}$. Moreover, since $YY'$ is positive semi-definite, $\tr((\hat{\Sigma}-\frac{1}{1-\delta})YY')\leq (\lambda_1-\frac{1}{1-\delta})\tr(YY')$. For $X\in D_n$, $(\lambda_1-\frac{1}{1-\delta})<\delta^{-1}p/n$. As a result,
\[
J_2 \leq e^{\frac{nr\delta}{2p(1-\delta)}} \int_{\delta\eta_1>1/2} \mexp\bigr\{\frac{1}{2}\tr(YY')\bigr\}dF(Y). 
\]
Let $W=(2pr)^{-1/2}(\tr(YY')-pr)$. Note that $\eta_1\geq (pr)^{-1}\tr(YY')$. So $W>\frac{1-\delta}{2\delta}\sqrt{(pr)/2}$ implies that $\delta\eta_1>1/2$. It follows that
\[
J_2 \leq e^{\frac{nr\delta}{2p(1-\delta)}+\frac{pr}{2}} \int_{W>\frac{1-\delta}{2\delta}\sqrt{\frac{pr}{2}}} e^{\sqrt{\frac{pr}{2}} W}dF(W). 
\]
Since $\tr(YY')$ has a chi-square distribution, $W$ weakly converges to $N(0,1)$. So $e^{\sqrt{\frac{pr}{2}}W}dF(W)=e^{\frac{pr}{2}}\cdot \frac{1}{\sqrt{2\pi}}e^{\frac{1}{2}(W-\sqrt{\frac{pr}{2}})^2}[1+o(1)]dW$, and the above integral is equal to $e^{pr/2}[1+o(1)] \cdot P\big(N(0,1)>(\frac{1-\delta}{2\delta}-1)\sqrt{pr/2}\big)$. It follows that $J_2\leq e^{\frac{nr\delta}{2p(1-\delta)}+\frac{pr}{2}}\cdot 2e^{ -Cpr/\delta^2}\leq 2e^{-Cpr/(2\delta^2)}=o(e_n)$. 

Consider $J_3$. First, $X\in D_n^c$ implies $\lambda_1>\frac{p}{n\delta}+\frac{1}{1-\delta}>\frac{p}{2n\delta}$, noting that $\delta\goto 0$. Second, under $\widetilde{H}_1^{(n)}$, $X_i|Y\overset{iid}{\sim} N(0, \widetilde{\Sigma}(Y))$, $1\leq i\leq n$, where $[\widetilde{\Sigma}(Y)]^{-1}=I_p-\delta p^{-1}YY'$. Since for $Y\in A_n$, the maximum eigenvalue of $\widetilde{\Sigma}(Y)$ is $(1-\delta\eta_1)^{-1}\leq 2$, we find that $\widetilde{P}_1(\lambda_1>t)\leq P_0(2\lambda_1>t)$ for any $t>0$, where $P_0$ is the probability measure under $H_0^{(n)}$. 
Combining the above results and applying \eqref{Vershynin5.35}, 
\[
J_3 \leq P_0\Big(2\lambda_1>\frac{p}{2n\delta}\Big)\leq 2e^{-\frac{n}{2}(\frac{p}{4n\delta}-1-\sqrt{\frac{p}{n}})^2}\leq 2 e^{- C\delta^{-1}n}=o(e_n). 
\]

Last, consider $J_1$. Let $q(x)=-\log(1-x)-x-x^2/2$ and recall that $\eta_1>\eta_2>\cdots>\eta_r>0$ are the nonzero eigenvalues of $(1/p)YY'$. Introduce
\[
M_1(Y) = \sum_{k=1}^r q\big(\frac{\delta(\eta_k-1)}{1-\delta}\big), \quad
M_2(Y) =\frac{\delta^2}{2(1-\delta)^2} \sum_{k=1}^r\bigl\{ (\eta_k-1)^2 - E[(\eta_k-1)^2]\bigr\}. 
\]
Then, $\Delta_n(Y)=M_1(Y)+M_2(Y)$. 
Let $E^*[f(Y)]=\int_{Y\in A_n}f(Y)dF(Y)$ for any function $f$. We then have
\begin{align}  \label{approx-J1}
J_1 &\leq E^*|1-e^{\frac{n}{2}M_1(Y)}| + E^*\bigl[e^{\frac{n}{2}M_1(Y)}|1-e^{\frac{n}{2}M_2(Y)}|\bigr]\cr
&\leq E^*|1-e^{\frac{n}{2}M_1(Y)}| + \bigl(E^*[ e^{nM_1(Y)}]\cdot E^*|1-e^{\frac{n}{2}M_2(Y)}|^2\bigr)^{1/2}. 
\end{align}
We claim that  
\beq \label{approx-5}
E^*|1-e^{\frac{n}{2}M_1(Y)}|\leq Ce_n, \qquad E^*|1-e^{nM_1(Y)}|\leq Ce_n,
\eeq
and 
\beq \label{approx-6}
E^*|1-e^{\frac{n}{2}M_2(Y)}|\leq Ce_n, \qquad E^*|1-e^{nM_2(Y)}|\leq Ce_n. 
\eeq
Note that \eqref{approx-6} implies $E^*|1-e^{\frac{n}{2}M_2(Y)}|^2\leq Ce_n$. 
Combining it with \eqref{approx-J1}-\eqref{approx-5} gives $J_1\leq Ce_n$.  

It remains to show \eqref{approx-5}-\eqref{approx-6}. We first prove \eqref{approx-5}. It suffices to consider the first inequality, as the second one is similar. Define $U=\frac{\delta}{1-\delta}(\eta_1-1)$ if $|\eta_1-1|\geq |\eta_r-1|$, and $U=\frac{\delta}{1-\delta}(\eta_r-1)$ if $|\eta_1-1|<|\eta_r-1|$. Note that $U\leq1/2$ for $Y\in A_n$. 
Since $q(x)=\sum_{\ell=3}^\infty x^\ell/\ell$, we have $|M_1(Y)|\leq rq(|U|)$. Therefore, to show \eqref{approx-5}, it suffices to show that 
\beq \label{approx-7}
E\bigl(|1-e^{\frac{nr}{2}q(|U|)}|\cdot 1\{U\leq 1/2\}\bigr)\leq Ce_n.  
\eeq
Let $b_n=3\delta\sqrt{r/p}$. By Taylor expansion, $q(|x|)\leq C_1|x|^3$ for $x\in[-b_n, b_n]$, $q(|x|)\leq C_2x^2$, for $x\in (b_n, 1/2]\cup (-b_n,\infty)$, where $C_1,C_2>0$ are absolute constants. 
Note that $nrb_n^3\leq Ce_n\goto 0$. Hence, when $|U|\leq b_n$, $|1-e^{\frac{nr}{2}q(|U|)}|\leq Cnrb_n^3$; when $|U|>b_n$, $|1-e^{\frac{nr}{2}q(|U|)}|\leq 1+e^{\frac{nr}{2}q(|U|)}\leq 1+e^{CnrU^2}$. It follows that the left hand side of \eqref{approx-7} is upper bounded by 
\beq  \label{approx-8}
Cnrb_n^3 + P(|U|>b_n) + E\bigr[ e^{CnrU^2}1\{|U|>b_n\}\bigr]. 
\eeq
The first term is $O(e_n)$. 
To bound the second term, by \eqref{Vershynin5.35}, 
write 
\[
|U|=\frac{\delta}{1-\delta}(\sqrt{r/p}+U_1/\sqrt{p}), \qquad \mbox{where }\;\;P(U_1>t)\leq 2e^{-t^2/2}\;\; \mbox{ for $t>0$}.
\]
It follows that $P(|U|>b_n)\leq P(U_1>\sqrt{r})\leq 2e^{-r^2/2}=o(e_n)$. 
Last, consider the third term. Since $nrU^2\leq C\delta^2r^2+C\delta^2rU_1^2$, this term is bounded by $e^{C\delta^2r^2}\cdot E[ e^{C\delta^2rU_1^2}1\{ U_1>\sqrt{r}\}]\leq e^{C\delta^2r^2-r/4}$, provided that $\sqrt{r}\gg \delta^2r$ (which is obviously true here). Since $\delta^2r\leq e_n\goto 0$, this term is $O(e^{-r/8})=o(e_n)$. Combing the above gives \eqref{approx-7}, and \eqref{approx-5} follows. 

Next, we prove \eqref{approx-6}. We only consider the first inequality, and the second one is similar. 
Let $U$ and $b_n$ be the same as above, and let $V=\sum_{k=1}^r\{ (\eta_k-1)^2 - E[(\eta_k-1)^2]\}$. It is seen that $n|M_2(Y)|= n\frac{\delta^2}{2(1-\delta)^2} |V|$. Moreover, since $\sum_{k=1}^r(\eta_k-1)^2\leq r\max\{|\eta_1-1|^2,|\eta_r-1|^2\}$ and $\sum_{k=1}^rE[(\eta_k-1)^2]=r/p$, we have a naive upper bound: $n|M_2(Y)|\leq CnrU^2+C\delta^2r$. As a result,
\begin{align}  \label{approx-9}
E^*|1-e^{\frac{n}{2}M_2(Y)}|&\leq E\bigl( |1-e^{Cn\delta^2 V}|\cdot 1\{|U|\leq b_n\}\bigr) + P(|U|>b_n)\cr
& + e^{C\delta^2r} E\bigl[ e^{CnrU^2} 1\{ |U|>b_n \}   \bigr]. 
\end{align}
In the analysis of \eqref{approx-8}, we have seen that the second term is $o(e_n)$ and the last term is $O(e^{C\delta^2r-r/8})=o(e_n)$. 
Consider the first term. 
We adapt the proof of Lemma~\ref{lem:eigen-moment} by constructing a function $\tilde{h}_n$, similar to $h_n$ in the proof of Lemma~\ref{lem:eigen-moment}, except that $a_n$ is now replaced with $3\sqrt{r/p}$. Then, $\|\tilde{h}_n'\|_\infty\leq 6\sqrt{r/p}$; in addition, when $|U|\leq b_n$, $V=\sum_{k=1}^r \tilde{h}_n(\eta_k)-\sum_{k=1}^r E[(\eta_k-1)^2]$. Similar to \eqref{lem-mmt-4}, we can show that $\sum_{k=1}^r E[(\eta_k-1)^2]=\sum_{k=1}^r E[\tilde{h}_n(\eta_k)] +o(r/p)$. Therefore, by Theorem~\ref{thm:deviation} (with $m=2$),  
\[
P\bigl( |V|>6p^{-1}r t , |U|\leq b_n\bigr) \leq 2\inf_{y\geq\sqrt{2\log(2)}}\bigl\{e^{-\frac{t^2}{c_0y^2}} + 2e^{-y^2}\bigr\}\leq 4e^{-t/\sqrt{c_0}}. 
\]
With $V_1=r^{-1}pV\cdot 1\{|U|\leq b_n\}$, the above implies $P(|V_1|>6t)\leq 4e^{-c_0^{-1/2}t}$, so $V_1$ has an exponential tail. The first term in \eqref{approx-9} is upper bounded by $E|1-e^{C\delta^2r V_1}|\leq C\delta^2r=O(e_n)$. Combining the above gives \eqref{approx-6}. 
\qed  

\subsection{Proof of Lemma~\ref{lem:tildeR}}
We first show \eqref{proxy1a}. Let $\hat{\Sigma}=\hat{Q}\Lambda \hat{Q}'$ be the eigen decomposition of $\hat{\Sigma}$, where $\Lambda=\mathrm{diag}(\lambda_1,\cdots,\lambda_p)$ and $\hat{Q}$ is an orthogonal matrix. Let $Z=\hat{Q}'Y$. Write $b^*_n(\delta,r)= (1-\delta)^{\frac{nr}{2}}b_n(\delta, r)$. We can rewrite 
\begin{align*}
\widetilde{R}_n &= b^*_n(\delta, r) \int \mexp\left\{\frac{\delta n}{2p}\tr\bigl(Z'(\Lambda-\frac{1}{1-\delta}I_n)Z\bigr)\right\}dF(Z)\cr
&=b^*_n(\delta, r)  \int \mexp\bigg\{\frac{\delta n}{2p}\sum_{k=1}^r\sum_{j=1}^p (\lambda_j-\frac{1}{1-\delta})Z^2(j,k) \bigg\}dF(Z).
\end{align*}
Conditioning on $X$,  $Z$ has $iid$ entries of $N(0,1)$. For $z\sim N(0,1)$ and $c<1$, $E[e^{cz^2/2}]=(1-c)^{-1/2}=e^{-\log(1-c)/2}$. Therefore, 
\[
\widetilde{R}_n = b^*_n(\delta, r) \mexp\Big\{-\frac{r}{2}\sum_{j=1}^p \log\bigl(1- \frac{\delta n}{p}(\lambda_j-\frac{1}{1-\delta})\bigr)\Big\}. 
\]

We then show \eqref{proxy1b}. By definition, we have $\sum_{j=1}^p \log(1- \frac{\delta n}{p}(\lambda_j-\frac{1}{1-\delta}))= \frac{n\delta}{p}\sum_{j=1}^p \psi_n(\lambda_j)=n\delta \int\psi_n(\lambda)\mu_n^X(d\lambda)$. Hence, given \eqref{proxy1a}, it suffices to show that 
\beq \label{lem-tildeR-1}
\int \psi_n(\lambda)\mu_{n,p}(d\lambda) = \frac{1}{\delta}\log(1-\delta) + \frac{1}{1-\delta}. 
\eeq
Write $\gamma_n=p/n$. Then $\psi_n(\lambda)=\frac{\gamma_n}{\delta}\log(\frac{\delta}{\gamma_n})+\frac{\gamma_n}{\delta}\log(\frac{\gamma_n}{\delta}+\frac{1}{1-\delta}-\lambda)$. 
Let $H_\gamma(x)=\int \log(x-\lambda)\mu_{\gamma}(d\lambda)$ for any $\gamma\geq 1$ and $x>(1+\sqrt{\gamma})^2$. It is seen that 
\beq  \label{lem-tildeR-2}
\int \psi_n(\lambda)\mu_{n,p}(d\lambda) = \frac{\gamma_n}{\delta} H_{\gamma_n}\bigl(\frac{\gamma_n}{\delta}+\frac{1}{1-\delta}\bigr) +\frac{\gamma_n}{\delta}\log(\delta/\gamma_n). 
\eeq
We apply Lemma~\ref{lem:MPintegral} with $t=\frac{\gamma_n}{\delta}+\frac{1}{1-\delta}$ and $\gamma=\gamma_n$. By direct calculations, $A_t=\frac{1-\delta}{\delta}\gamma_n + \frac{\delta}{1-\delta}$, $B_t= \frac{1-\delta}{\delta}\gamma_n - \frac{\delta}{1-\delta}$ and $(\gamma_n+1)A_t-(\gamma_n-1)B_t+4\gamma_n = \frac{2\gamma_n}{\delta(1-\delta)}$. Plugging them into Lemma~\ref{lem:MPintegral} and inserting the result into \eqref{lem-tildeR-2} gives \eqref{lem-tildeR-1}. 
\qed

\section{Acknowlegement} 
The author thanks Jiashun Jin for pointing out the topic and numerous thoughtful pointers and comments on the manuscript. She thanks Zhigang Bao for personal communications on Random Matrix Theory. She thanks Mark Low, Zhonggen Su and Ofer Zeitouni for helpful pointers.   


\bibliographystyle{imsart-number}
\bibliography{Jiashun}
\end{document}